 \documentclass[a4paper,12pt]{article}
    \usepackage[top=2.5cm,bottom=2.5cm,left=2.5cm,right=2.5cm]{geometry}
    \usepackage{cite, amsmath, amssymb}
   \usepackage{tikz}
\begin{document}
\begin{center}
{\LARGE\bf Integral cographs and applications}
\end{center}
\begin{center}
{\large \bf Luiz Emilio Allem$^a,$ Fernando Tura$^{b,}$ \footnote{Corresponding author} }
\end{center}
\begin{center}
\it $^a$ Instituto  de Matem\'atica, UFRGS, Porto Alegre, RS, 91509-900, Brazil\\
\tt emilio.allem@ufrgs.br
\end{center}
\begin{center}
\it $^b$ Departamento de Matem\'atica, UFSM, Santa Maria, RS, 97105-900, Brazil\\
\tt ftura@smail.ufsm.br
\end{center}

\newcommand{\lambdamin}{\lambda_{\min, n}}
\newcommand{\formulamin}{
 \frac{  (
   \lfloor \frac{n}{3} \rfloor \! - \! 1 ) \! - \! \sqrt{ ( \lfloor \frac{n}{3} \rfloor \! - \! 1)^2 \! +\!  4
   (n \! - \! \lfloor \frac{n}{3} \rfloor )
   \lfloor \frac{n}{3} \rfloor
  }
  }{2}
}
\newcommand{\casei}{{\bf case~1}}
\newcommand{\subia}{{\bf subcase~1a}}
\newcommand{\subib}{{\bf subcase~1b}}
\newcommand{\subic}{{\bf subcase~1c}}
\newcommand{\caseii}{{\bf case~2}}
\newcommand{\subiia}{{\bf subcase~2a}}
\newcommand{\subiib}{{\bf subcase~2b}}
\newcommand{\caseiii}{{\bf case~3}}
\newcommand{\myvar}{x}
\newcommand{\exvar}{\frac{\sqrt{3} + 1}{2}}
\newcommand{\PrfSketch}{{\bf Proof (Sketch): }}
\newcommand{\boldQ}{\mbox{\bf Q}}
\newcommand{\boldR}{\mbox{\bf R}}
\newcommand{\boldZ}{\mbox{\bf Z}}
\newcommand{\boldc}{\mbox{\bf c}}
\newcommand{\sign}{\mbox{sign}}
\newcommand{\alphaseq}{{\pmb \alpha}_{G,\myvar}}
\newcommand{\alphaseqlambda}{{\pmb \alpha}_{G,\lambda}}
\newcommand{\alphaseqGprime}{{\pmb \alpha}_{G^\prime,\myvar}}
\newcommand{\alphaseqlam}{{\pmb \alpha}_{G,-\lambdamin}}
\newtheorem{Thr}{Theorem}
\newtheorem{Pro}{Proposition}
\newtheorem{Que}{Question}
\newtheorem{Con}{Conjecture}
\newtheorem{Cor}{Corollary}
\newtheorem{Lem}{Lemma}
\newtheorem{Fac}{Fact}
\newtheorem{Ex}{Example}
\newtheorem{Def}{Definition}
\newtheorem{Prop}{Proposition}
\newtheorem{Remark}{Remark}

\def\floor#1{\left\lfloor{#1}\right\rfloor}

\newenvironment{my_enumerate}{
\begin{enumerate}
  \setlength{\baselineskip}{14pt}
  \setlength{\parskip}{0pt}
  \setlength{\parsep}{0pt}}{\end{enumerate}
}
\newenvironment{my_description}{
\begin{description}
  \setlength{\baselineskip}{14pt}
  \setlength{\parskip}{0pt}
  \setlength{\parsep}{0pt}}{\end{description}
}

\begin{abstract}
A graph is called integral if all the eigenvalues of its adjacency matrix are  integers. In this paper, we show a cograph that has a balanced cotree $T_{G}(a_{1},\ldots,a_{r-1},0|0,\ldots,0,a_{r})$ is integral computing its spectrum. As an application, these integral cographs can be used to estimate the eigenvalues of any cograph.
\\

\noindent
{\bf keywords:} cograph, adjacency  matrix, integral graphs, interlacing theorem.  
\end{abstract}



\section{Introduction}
\label{intro}

 The spectrum of a graph $G$ having $n$ vertices is the (multi)set of the eigenvalues of its adjacency matrix. We traditionally order them so that
$$ \lambda_1 \geq \lambda_2 \geq \cdots \geq \lambda_n.$$
A graph is called integral if all the eigenvalues are  integers. The notion of integral graphs dates back to Harary  and Schwenk \cite{Harary1}. Since then, se\-ve\-ral explicit constructions of integral graphs of special types  appears in the literature, see for example \cite{Cve1, Bro1,Bro2,Moha1}. In \cite{bali}, the authors present, using computational experiments, that out of $164,059,830,476$
connected graphs on $12$ vertices, there exist exactly $325$ integral graphs. It shows finding integral graphs is a very difficult task. And, an interesting fact is that, integral graphs have applications in quantum networks allowing perfect state transfer \cite{Saxena}.

 In this paper we deal with cographs, an important class of graphs for its many applications. There are many ways to characterize cographs. For example,  a cograph is a graph which contains no  path  of length  four  as an induced subgraph \cite{Stewart}, that is why they are often simply  called  $P_4$ free graph. Any cograph has  a canonical tree representation, called  the  cotree \cite{BSS2011}.

There is a considerable body of knowledge on the spectral properties of graphs \cite{Bro}. As far as cographs, spectral properties were studied by  \cite{BSS2011}. They determined the multiplicity of the eigenvalues  $-1$ and $0$. Royle in \cite{Royle} proved  that the rank of the adjacency matrix of any cograph is equal to the number of distinct non-zero rows of that.
Spectral properties of threshold graphs (a subclass of cographs) were studied in \cite{JTT2015, JTT2013, tura18}.
In  \cite{Moha} was proved  that no cograph has eingenvalues in the interval $(-1,0)$, a surprising result.

Although there is a lot of constructions of integral graphs, the literature does not seem to provide  many results about integral cographs.
In this paper, we attempt to fill this gap, presenting families of integral cographs.

The main tool used here is an algorithm called Diagonalize, presented in \cite{JTT2016}.  The algorithm
finds, in $O(n)$ time, the  number of eigenvalues that are greater than $x$, less than $x$ and equal to $x$, operating directly on the cotree of the cograph.

Here is an outline of this paper.
In section \ref{Diag}, Algorithm Diagonalize and some of its properties are explained.
In section \ref{balancedcotree}, some characteristics of balanced cotrees are presented.
In section \ref{recurrence}, we compute eigenvalues of the cograph that has balanced cotree using recurrence relations.
In section \ref{multi}, we give the multiplicity of the eigenvalues computed in section \ref{recurrence}.
In section \ref{examplesection}, examples of integral cographs are presented.
In section \ref{applic}, we finalize this paper presenting how to estimate the eigenvalues of any cograph by Interlacing Theorem using the integral cographs with balanced cotree $T_{G}(a_{1},\ldots,a_{r-1},0|0,\ldots,0,a_{r})$.


\section{Preliminaries}
\label{Diag}
Let  $G= (V,E)$ be an undirected graph with vertex set $V$ and edge set $E$, without loops or multiple edges.  For $v\in V$, $N(v)$ denotes the  open neighborhood  of $v$, that is, $\{w|\{v,w\}\in E\}$. The  closed neighborhood  $N[v] = N(v) \cup \{v\}$. If $|V| = n$, the   adjacency matrix  $A= (a_{ij})$,  is
 the $n \times n$   matrix of zeros and ones  such that $a_{ij} = 1$ if
there is an edge between $v_i$ and $v_j$, and 0 otherwise.  A value $\lambda$ is an
{\em eigenvalue} of $G$ if $\det(A - \lambda I_n ) = 0$, and since
$A$ is a matrix real and symmetric, its eigenvalues are real numbers.

A cograph has been rediscovered independently by  several authors
since the 1960's.  Corneil, Lerchs and Burlingham \cite{Stewart} define cographs recursively by the following rules:
\begin{enumerate}
	\item [(i)] a graph on a single vertex is a cograph,
  \item [(ii)] a finite  union and join of cographs are a cograph.
\end{enumerate}

\begin{figure}[h!]
       \begin{minipage}[c]{0.3 \linewidth}
\begin{tikzpicture}
  [scale=0.6,auto=left,every node/.style={circle}]
  \foreach \i/\w in {1/,2/,3/,4/,5/,6/,7/}{
    \node[draw,circle,fill=black,label={360/7 * (\i - 1)+90}:\i] (\i) at ({360/7 * (\i - 1)+90}:3) {\w};} 
  \foreach \from in {4}{
    \foreach \to in {3,\from}
      \draw (\from) -- (\to);}
       \foreach \from in {2}{
       \foreach \to in {1,\from}
       \draw (\from) -- (\to)
       ;}

\foreach \from in {5}{
       \foreach \to in {1,2,3,4,\from}
       \draw (\from) -- (\to)
       ;}

\foreach \from in {6}{
       \foreach \to in {1,2,3,4,5,\from}
       \draw (\from) -- (\to)
       ;}

\foreach \from in {7}{
       \foreach \to in {1,2,3,4,\from}
       \draw (\from) -- (\to)
       ;}

\end{tikzpicture}
       \end{minipage}\hfill
       \begin{minipage}[l]{0.4 \linewidth}
\begin{tikzpicture}
 [scale=1,auto=left,every node/.style={circle,scale=0.9}]

  \node[draw,circle,fill=black,label=below:$v_6$] (o) at (3,4) {};
  \node[draw,circle,fill=black,label=below:$v_4$] (n) at (1.8,4) {};
   \node[draw,circle,fill=black,label=below:$v_5$] (q) at (2.2,4) {};

  \node[draw, circle, fill=blue!10, inner sep=0] (m) at (2.5,5) {$\otimes$};
  \node[draw,circle,fill=black,label=below:$v_7$] (l) at (4,5) {};

  \node[draw, circle, fill=blue!10, inner sep=0] (j) at (3,6) {$\cup$};
  \node[draw,circle,fill=blue!10, inner sep=0] (h) at (2,7) {$\otimes$};
  \node[draw, circle, fill=blue!10, inner sep=0] (g) at (1,6) {$\cup$};
  \node[draw,circle,fill=blue!10, inner sep=0] (f) at  (1.3,5) {$\otimes$};
  \node[draw,circle,fill=blue!10, inner sep=0] (a) at (0,5) {$\otimes$};
  \node[draw, circle, fill=black, label=below:$v_1$] (b) at (-0.3,4) {};
  \node[draw,circle,fill=black, label=below:$v_2$] (c) at (0.5,4) {};
  \node[draw,circle,fill=black,label=below:$v_3$] (e) at (1,4) {};

  \path (a) edge node[left]{} (b)
        (a) edge node[below]{} (c)

   (f) edge node[below]{} (e)
  (f) edge node[below]{} (n)


        (f) edge node[right]{}(g)
        (g) edge node[left]{}(a)
        (h) edge node[right]{}(j)
        (h) edge node[left]{}(g)

        (j) edge node[right]{}(l)
        (j) edge node[below]{}(m)
        (m) edge node[right]{}(o)
        (m) edge node[left]{} (q);
\end{tikzpicture}
       \end{minipage}
       \caption{A cograph $G=((v_{1}\otimes v_{2})\cup (v_{3}\otimes v_{4}))\otimes ((v_{5}\otimes v_{6}) \cup v_7))$ and its cotree $T_G$.}
       \label{cotree}
\end{figure}
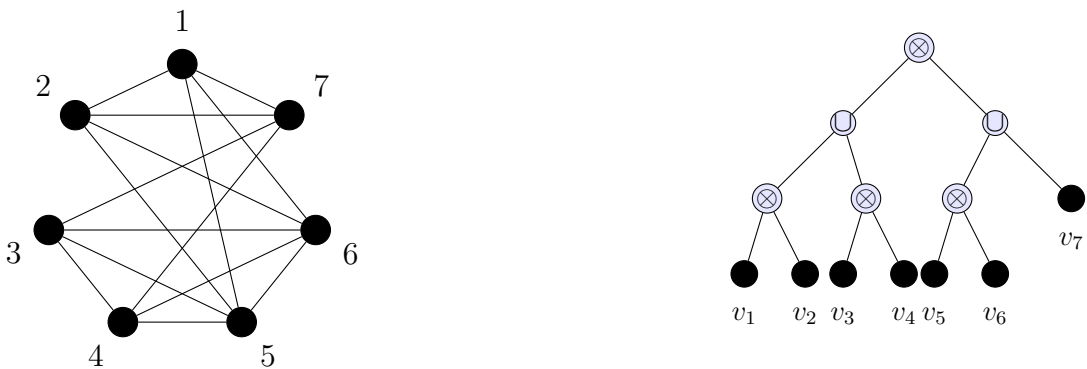

A cotree $T_G$ of a cograph $G$ is a rooted tree in which any interior vertex $w$ is either of $\cup$ type (corresponds to union) or $\otimes$ type (corresponds to join). The terminal vertices (leaves) are typeless and represent the vertices of the cograph $G$. We say that {\em depth} of the cotree is the number of edges of  the longest path from the root to a leaf. To build a cotree for a connected cograph, we simply place a $\otimes$ at the tree's root, placing $\cup$ on interior vertices with odd depth, and placing $\otimes$ on interior vertices with even depth. All interior vertices have at least two children.
The Figure \ref{cotree} shows  a cograph and its cotree.

Two vertices $u$ and $v$  are duplicate if $N(u) = N(v)$ and coduplicate if $N[u]=N[v]$.
An important tool presented in \cite{JTT2016} was an algorithm
for constructing a {\em diagonal} matrix congruent to $A + \myvar I$,
where $A$ is the adjacency matrix of a cograph,
and $\myvar$ is an arbitrary scalar, using $O(n)$ time and space.

The algorithm's input  is the cotree $T_G$ and $x$. Each leaf $v_i$, $i =1,\ldots,n$ have a value $d_i$ that represents the diagonal element of  $A+xI$. It initializes  all entries $d_i$ with $x$.
At each iteration, a pair $\{v_k, v_l\} $ of duplicate or coduplicate vertices with maximum depth  is selected. Then they are processed, that is,  assignments are given to $d_k$ and $d_l$, such that either one or both rows (columns) are diagonalized. When a $k$ row (column)  corresponding to vertex  $v_k$  has been diagonalized then $v_k$ is  removed from the $T_G$, it means that $d_k$ has a permanent final value. Then the algorithm moves to the cotree $T_G -v_k$. The algorithm is shown in Figure \ref{algo}.

\begin{figure}[h]
{\small
{\tt
\begin{tabbing}
aaa\=aaa\=aaa\=aaa\=aaa\=aaa\=aaa\=aaa\= \kill
     \> INPUT:  cotree $T_G$, scalar $\myvar$\\
     \> OUTPUT: diagonal matrix $D=[d_1, d_2, \ldots, d_n]$ congruent to $A(G) + \myvar I$\\
     \>\\
     \>   $\mbox{ Algorithm}$ Diagonalize $(T_{G}, x)$ \\
     \> \> initialize $d_i := \myvar$, for $ 1 \leq i \leq n$ \\
     \> \> {\bf while } $T_G$  has $\geq 2$    leaves      \\
     \> \> \>  select a pair $(v_k, v_l)$  (co)duplicate of maximum depth with  parent $w$\\
     \> \> \>     $\alpha \leftarrow  d_k$    $\beta \leftarrow d_{l}$\\
     \> \> \> {\bf if} $ w=\otimes$\\
     \> \> \> \>  {\bf if} $\alpha + \beta \neq2$  \verb+                //subcase 1a+    \\
     \> \> \> \> \>   $d_{l} \leftarrow \frac{\alpha \beta -1}{\alpha + \beta -2};$ \hspace*{0,25cm} $d_{k} \leftarrow \alpha + \beta -2; $\hspace{0,25cm}   $T_G = T_G - v_k$ \\
     \> \> \> \>  {\bf else if } $\beta=1$ \verb+                //subcase 1b+   \\
     \> \> \> \> \>   $d_{l} \leftarrow 1$ \hspace*{0,25cm}   $d_k  \leftarrow 0;$ \hspace{0,25cm} $T_G = T_G - v_k$ \\
     \> \> \> \>  {\bf else  }  \verb+                      //subcase 1c+   \\
     \> \> \> \> \>   $d_{l} \leftarrow 1$  \hspace*{0,25cm} $d_k \leftarrow -(1-\beta)^2;$ \hspace{0,25cm} $T_G= T_G -v_k;$ \hspace{0,25cm} $T_G = T_G -v_l$  \\
     \> \> \>     {\bf else if} $w=\cup$\\
     \> \> \> \>  {\bf if} $\alpha + \beta \neq 0$  \verb+               //subcase 2a+    \\
     \> \> \> \> \>   $d_{l} \leftarrow \frac{\alpha \beta}{\alpha +\beta};$ \hspace*{0,25cm}   $d_k \leftarrow \alpha +\beta;$ \hspace{0,25cm} $T_G = T_G - v_k$ \\
     \> \> \> \>  {\bf else if } $\beta=0$ \verb+                //subcase 2b+   \\
     \> \> \> \> \>   $d_{l} \leftarrow 0;$ \hspace*{0,25cm}  $d_k  \leftarrow 0;$ \hspace{0,25cm} $T_G = T_G - v_k$ \\
     \> \> \> \>  {\bf else  }  \verb+                      //subcase 2c+   \\
     \> \> \> \> \>   $d_{l} \leftarrow \beta;$  \hspace*{0,25cm} $v_k \leftarrow -\beta;$ \hspace{0,25cm} $T_G =T_G - v_k;$ \hspace{0,25cm} $T_G = T_G - v_l$  \\

     \> \>  {\bf end loop}\\
\end{tabbing}
}}
\caption{\label{algo} Algorithm Diagonalize.}
\end{figure}

The next three results from \cite{JTT2016} will be used throughout the paper.

\begin{Thr}
\label{main1}
Let $D=[d_1,d_2,\ldots,d_n]$ be the diagonal returned  by the Algorithm Diagonalize $(T_G,-x)$, and assume $D$ has $k_{+}$ positive values, $k_0$ zeros and $k_{-}$ negative  values.
\begin{my_description}
 \item[i]
The number of eigenvalues of $G$ that are greater than $\myvar$ is exactly $k_{+}$.
\item[ii]
The number of eigenvalues of $G$ that are less than $\myvar$ is exactly $k_{-}$.
\item[iii]
The multiplicity of  $\myvar$ is  $k_{0}$.
\end{my_description}\end{Thr}

The following two lemmas show if a node $\otimes$ or $\cup$, in the cotree, have leaves with the same value. Then, the assignments can be controlled directly.

\begin{Lem}
\label{lem1}
If $v_1, \ldots, v_m$ have parent $w= \otimes$, each with diagonal value $y \neq 1$, then the algorithm performs $m-1$ iterations of  \subia~ assigning, during iteration  $j:$
\begin{equation}
d_k  \leftarrow \frac{j+1}{j}(y-1)
\end{equation}
\begin{equation}
d_l  \leftarrow \frac{y+j}{j+1}
\end{equation}
\end{Lem}

\begin{Lem}
\label{lem2}
If $v_1, \ldots, v_m$ have parent $w= \cup$, each with diagonal value $y\neq 0$, then the algorithm performs $m-1$ iterations of  \subiia~ assigning, during iteration  $j:$
\begin{equation}
d_k  \leftarrow \frac{(j+1)}{j}y
\end{equation}
\begin{equation}
d_l  \leftarrow \frac{y}{j+1}
\end{equation}
\end{Lem}


\section{Balanced Cotree $T_{G}(a_{1},\ldots,a_{r-1},0|0,\ldots,0,a_{r})$.}
\label{balancedcotree}

The balanced cotree $T_{G}(a_{1},\ldots,a_{r-1},0|0,\ldots,0,a_{r})$ of depth $r$ has a node $\otimes$ at the root, this node has exactly $a_1$ immediate $\cup$ interior vertices. Each $\cup$ at level $1$ has exactly $a_2$ immediate $\otimes$ interior vertices, and so on. Notice that, this cotree only has leaves at the last level. It means that, the nodes at level $r-1$ have $a_{r}$ immediate leaves. So, at level $i$, the cotree has $a_{1}a_{2}\cdots a_{i}$ nodes $\otimes$ if $i$ is even and $\cup$ if $i$ is odd, for $1\leq i\leq r-1$. And, at the last level $r$, it has $a_{1}a_{2}\cdots a_{r}$ leaves.

\begin{figure}[h!]
\begin{center}
\begin{tikzpicture}
  [scale=1,auto=left,every node/.style={circle,scale=0.9}]

        \node[draw, circle, fill=blue!10, inner sep=0] (m5) at (7.7,8) {$\otimes$};
  \node[draw,circle,fill=blue!10, inner sep=0] (l5) at (8.4,8) {$\otimes$};
  \node[draw,circle,fill=blue!10, inner sep=0] (y5) at (9.5,9) {$\cup$};
   \node[draw,circle,fill=blue!10, inner sep=0] (y25) at (9.5,8) {$\otimes$};
  \node[draw,circle,fill=blue!10, inner sep=0] (k5) at (10.5,8) {$\otimes$};
  \node[draw, circle, fill=blue!10, inner sep=0] (j5) at (8,9) {$\cup$};
  \node[draw,circle,fill=blue!10, inner sep=0] (h5) at (8,10) {$\otimes$};
  \node[draw, circle, fill=blue!10, inner sep=0] (g5) at (6.5,9) {$\cup$};
  \node[draw,circle,fill=blue!10, inner sep=0] (f5) at  (7,8) {$\otimes$};
  \node[draw,circle,fill=blue!10, inner sep=0] (a5) at (6,8) {$\otimes$};

  \node[draw,circle,fill=black, label=below:$$] (a6) at (5.75,7) {};
  \node[draw,circle,fill=black, label=below:$$] (b6) at (6.25,7) {};
  \node[draw,circle,fill=black, label=below:$$] (c6) at (6.75,7) {};
  \node[draw,circle,fill=black, label=below:$$] (d6) at (7.25,7) {};
   \node[draw,circle,fill=black, label=below:$$] (e6) at (7.65,7) {};
    \node[draw,circle,fill=black, label=below:$$] (f6) at (8.1,7) {};
     \node[draw,circle,fill=black, label=below:$$] (g6) at (8.5,7) {};
 \node[draw,circle,fill=black, label=below:$$] (h6) at (9,7) {};
  \node[draw,circle,fill=black, label=below:$$] (i6) at (9.5,7) {};
   \node[draw,circle,fill=black,label=below:$$] (j6) at (10,7) {};
    \node[draw,circle,fill=black, label=below:$$] (k6) at (10.5,7) {};
  \node[draw,circle,fill=black, label=below:$$] (l6) at (11,7) {};
 \path

                    (a5) edge node[right]{}(a6)
                    (a5) edge node[right]{}(b6)
      (f5) edge node[right]{}(c6)
      (f5) edge node[right]{}(d6)
               (m5) edge node[right]{}(e6)
      (m5) edge node[right]{}(f6)
                        (l5) edge node[right]{}(g6)
                        (l5) edge node[right]{}(h6)

                        (y25) edge node[right]{}(i6)
                        (y25) edge node[right]{}(j6)
                        (k5) edge node[right]{}(k6)
                        (k5) edge node[right]{}(l6)

             (h5) edge node[right]{}(y5)
        (k5) edge node[right]{}(y5)
        (y25) edge node[right]{}(y5)

        (f5) edge node[right]{}(g5)
        (g5) edge node[left]{}(a5)

        (h5) edge node[left]{}(j5)
        (h5) edge node[left]{}(g5)

        (j5) edge node[right]{}(l5)
        (j5) edge node[below]{}(m5)
        ;
\end{tikzpicture}
\caption{ $T_G(3,2,0|0,0,2)$.}
       \label{fig3}
\end{center}
\end{figure}
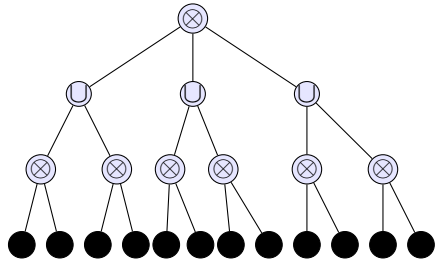

For more details on balanced cotrees see \cite{Allem}. In Figure \ref{fig3}, the balanced cotree $T_G (3,2,0| 0,0,2)$ with depth $r=3$ is shown. The complete graph $K_n$ is a cograph with balanced cotree $T_G(0|n)$.

The multiplicity of the eigenvalue $\lambda$ of the graph $G$ is denoted by $m(\lambda,G)$.
The next theorem is known and can be found in \cite{Allem}.

\begin{Thr}
\label{main2}
Let $G$ be a cograph with balanced cotree $T_{G}$ $(a_1,  \ldots,a_{r-1}, 0|$ $ 0, \ldots, 0,a_r)$ of order $n =a_{1} a_2  \ldots  a_{r-1} a_r$.
\begin{enumerate}
 \item[(i)]
 If $r$ is odd then $G$ has the eigenvalue $-1$ with
 \begin{equation}
 m(-1,G)= a_1 a_2 \ldots a_{r-1}(a_{r} -1)
 \end{equation}
\item[(ii)]
 If $r$ is even then $G$ has the eigenvalue $0$ with
 \begin{equation}
 m(0,G)=a_1 a_2 \ldots a_{r-1}(a_{r} -1)
 \end{equation}
\end{enumerate}
\end{Thr}

 For example, let $T_{G}$ be the balanced cotree in Figure \ref{fig3}. Then, $m(-1,G)=a_{1}a_{2}(a_{3}-1)=(3)(2)(2-1)=6$ and $m(0,G)=0$ by Theorem \ref{main2}.

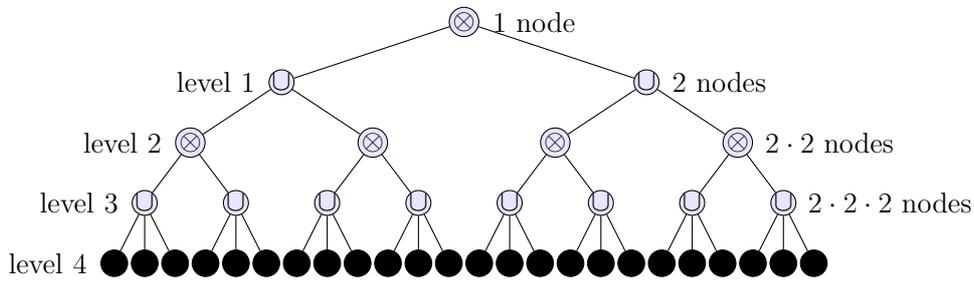
\begin{figure}[h!]
\begin{center}
\begin{tikzpicture}
  [scale=0.8,auto=left,every node/.style={circle,scale=0.9}]

                \node[draw, circle, fill=blue!10, inner sep=0,label=right:$1$ node] (a) at (5.75,4) {$\otimes$};
								
              \node[draw, circle, fill=blue!10, inner sep=0,label=left:level 1] (b) at (2.75,3) {$\cup$};
							 \node[draw, circle, fill=blue!10, inner sep=0,label=right:$2$ nodes] (c) at (8.75,3) {$\cup$};

        \node[draw, circle, fill=blue!10, inner sep=0,,label=left:level 2] (d) at (1.25,2) {$\otimes$};
				\node[draw, circle, fill=blue!10, inner sep=0] (e) at (4.25,2) {$\otimes$};
				\node[draw, circle, fill=blue!10, inner sep=0] (f) at (7.25,2) {$\otimes$};
				\node[draw, circle, fill=blue!10, inner sep=0,label=right:$2\cdot 2$ nodes] (g) at (10.25,2) {$\otimes$};

       \node[draw, circle, fill=blue!10, inner sep=0,label=left:level 3] (h) at (0.5,1) {$\cup$};
       \node[draw, circle, fill=blue!10, inner sep=0] (i) at (2,1) {$\cup$};			
			\node[draw, circle, fill=blue!10, inner sep=0] (j) at (3.5,1) {$\cup$};
			\node[draw, circle, fill=blue!10, inner sep=0] (k) at (5,1) {$\cup$};
			\node[draw, circle, fill=blue!10, inner sep=0] (l) at (6.5,1) {$\cup$};
			\node[draw, circle, fill=blue!10, inner sep=0] (m) at (8,1) {$\cup$};
			\node[draw, circle, fill=blue!10, inner sep=0] (n) at (9.5,1) {$\cup$};
			\node[draw, circle, fill=blue!10, inner sep=0,label=right:$2\cdot 2\cdot 2$ nodes] (o) at (11,1) {$\cup$};
 			
			\node[draw,circle,fill=black,,label=left:level 4] (h1) at (0,0) {};
			\node[draw,circle,fill=black] (h2) at (0.5,0) {};
			\node[draw,circle,fill=black] (h3) at (1,0) {};
			
			\node[draw,circle,fill=black] (i1) at (1.5,0) {};
			\node[draw,circle,fill=black] (i2) at (2,0) {};
			\node[draw,circle,fill=black] (i3) at (2.5,0) {};
			
			\node[draw,circle,fill=black] (j1) at (3,0) {};
			\node[draw,circle,fill=black] (j2) at (3.5,0) {};
			\node[draw,circle,fill=black] (j3) at (4,0) {};
			
			\node[draw,circle,fill=black] (k1) at (4.5,0) {};
			\node[draw,circle,fill=black] (k2) at (5,0) {};
			\node[draw,circle,fill=black] (k3) at (5.5,0) {};
			
			\node[draw,circle,fill=black] (l1) at (6,0) {};
			\node[draw,circle,fill=black] (l2) at (6.5,0) {};
			\node[draw,circle,fill=black] (l3) at (7,0) {};
			
			\node[draw,circle,fill=black] (m1) at (7.5,0) {};
			\node[draw,circle,fill=black] (m2) at (8,0) {};
			\node[draw,circle,fill=black] (m3) at (8.5,0) {};
			
			\node[draw,circle,fill=black] (n1) at (9,0) {};
			\node[draw,circle,fill=black] (n2) at (9.5,0) {};
			\node[draw,circle,fill=black] (n3) at (10,0) {};
			
			\node[draw,circle,fill=black] (o1) at (10.5,0) {};
			\node[draw,circle,fill=black] (o2) at (11,0) {};
			\node[draw,circle,fill=black] (o3) at (11.5,0) {};

			  \path
				  (a) edge node[right]{}(b)
					(a) edge node[right]{}(c)
					
					(b) edge node[right]{}(d)
					(b) edge node[right]{}(e)
					
					(c) edge node[right]{}(f)
					(c) edge node[right]{}(g)
				
				  (d) edge node[right]{}(h)
					(d) edge node[right]{}(i)
					
					(e) edge node[right]{}(j)
					(e) edge node[right]{}(k)
					
					(f) edge node[right]{}(l)
					(f) edge node[right]{}(m)
					
					(g) edge node[right]{}(n)
					(g) edge node[right]{}(o)
				
				  (h) edge node[right]{}(h1)
					(h) edge node[right]{}(h2)
					(h) edge node[right]{}(h3)
					
          (i) edge node[right]{}(i1)
					(i) edge node[right]{}(i2)
					(i) edge node[right]{}(i3)
					
					(j) edge node[right]{}(j1)
					(j) edge node[right]{}(j2)
					(j) edge node[right]{}(j3)
					
					(k) edge node[right]{}(k1)
					(k) edge node[right]{}(k2)
					(k) edge node[right]{}(k3)
					
					(l) edge node[right]{}(l1)
					(l) edge node[right]{}(l2)
					(l) edge node[right]{}(l3)
					
					(m) edge node[right]{}(m1)
					(m) edge node[right]{}(m2)
					(m) edge node[right]{}(m3)
					
					(n) edge node[right]{}(n1)
					(n) edge node[right]{}(n2)
					(n) edge node[right]{}(n3)
					
					(o) edge node[right]{}(o1)
					(o) edge node[right]{}(o2)
					(o) edge node[right]{}(o3)
							
        ;
\end{tikzpicture}
\caption{$T_{G_{4}}(2,2,2,0| 0,0,0,3)$ with $2\cdot 2\cdot 2\cdot 3$ leaves.}
       \label{cotree_example1}
       \end{center}
	\end{figure}

From now on, we will be dealing with cographs $G$ that have balanced cotrees of type $T_{G}(a_{1},\ldots,a_{r-1},0|0,\ldots,0,a_{r})$, with $a_{i}\geq 2$ for $1\leq i\leq r$. The cotree $T_{G}(a_{1},\ldots,a_{r-1},0|0,\ldots,0,a_{r})$ will be denoted with subindex to represent the depth $r$. That is,
$$T_{G}(a_{1},\ldots,a_{r-1},0|0,\ldots,0,a_{r})=T_{G_{r}}(a_{1},\ldots,a_{r-1},0|0,\ldots,0,a_{r}),$$
or $T_{G_{r}}$ for short.
In Figure \ref{cotree_example1}, the balanced cotree $T_{G_{4}}(2,2,2,0|$ $0,0,0,3)$ of the cograph $G$ of order $n =2\cdot 2\cdot 2\cdot 3=24$ is shown.

The next example is very important because it explains the notation used in the article.

\begin{Ex}
We  apply Algorithm Diagonalize to the cotree in Figure \ref{cotree_example1} with $x=-3$. That is, the input is $(T_{G_{4}}(2,2,2,0| 0,0,0,3),-3)$. The algorithm works from the bottom up and initially each leave at level 4 receives $x_{4}=-3$, as shown in Figure \ref{cotree_example2}.

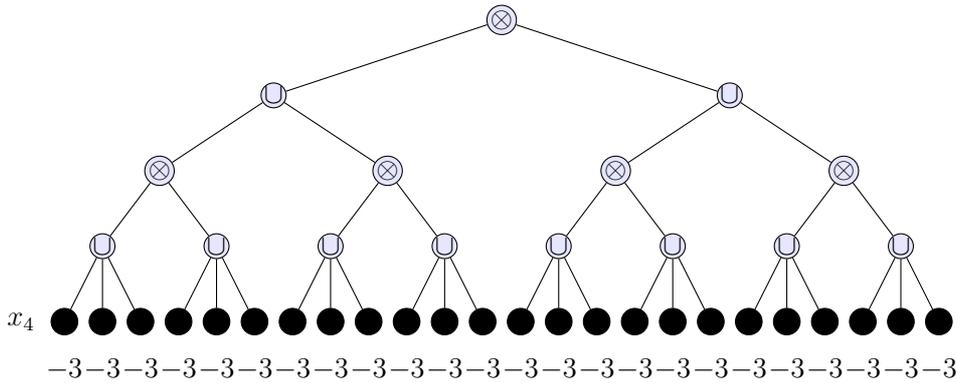
\begin{figure}[h!]
\begin{center}
\begin{tikzpicture}
  [scale=1,auto=left,every node/.style={circle,scale=0.9}]

                \node[draw, circle, fill=blue!10, inner sep=0] (a) at (5.75,4) {$\otimes$};
								
              \node[draw, circle, fill=blue!10, inner sep=0] (b) at (2.75,3) {$\cup$};
							 \node[draw, circle, fill=blue!10, inner sep=0] (c) at (8.75,3) {$\cup$};

        \node[draw, circle, fill=blue!10, inner sep=0] (d) at (1.25,2) {$\otimes$};
				\node[draw, circle, fill=blue!10, inner sep=0] (e) at (4.25,2) {$\otimes$};
				\node[draw, circle, fill=blue!10, inner sep=0] (f) at (7.25,2) {$\otimes$};
				\node[draw, circle, fill=blue!10, inner sep=0] (g) at (10.25,2) {$\otimes$};

       \node[draw, circle, fill=blue!10, inner sep=0] (h) at (0.5,1) {$\cup$};
       \node[draw, circle, fill=blue!10, inner sep=0] (i) at (2,1) {$\cup$};			
			\node[draw, circle, fill=blue!10, inner sep=0] (j) at (3.5,1) {$\cup$};
			\node[draw, circle, fill=blue!10, inner sep=0] (k) at (5,1) {$\cup$};
			\node[draw, circle, fill=blue!10, inner sep=0] (l) at (6.5,1) {$\cup$};
			\node[draw, circle, fill=blue!10, inner sep=0] (m) at (8,1) {$\cup$};
			\node[draw, circle, fill=blue!10, inner sep=0] (n) at (9.5,1) {$\cup$};
			\node[draw, circle, fill=blue!10, inner sep=0] (o) at (11,1) {$\cup$};
 			
			\node[draw,circle,fill=black,label=left:$x_{4}$,label=below:$-3$] (h1) at (0,0) {};
			\node[draw,circle,fill=black,label=below:$-3$] (h2) at (0.5,0) {};
			\node[draw,circle,fill=black,label=below:$-3$] (h3) at (1,0) {};
			
			\node[draw,circle,fill=black,label=below:$-3$] (i1) at (1.5,0) {};
			\node[draw,circle,fill=black,label=below:$-3$] (i2) at (2,0) {};
			\node[draw,circle,fill=black,label=below:$-3$] (i3) at (2.5,0) {};
			
			\node[draw,circle,fill=black,label=below:$-3$] (j1) at (3,0) {};
			\node[draw,circle,fill=black,label=below:$-3$] (j2) at (3.5,0) {};
			\node[draw,circle,fill=black,label=below:$-3$] (j3) at (4,0) {};
			
			\node[draw,circle,fill=black,label=below:$-3$] (k1) at (4.5,0) {};
			\node[draw,circle,fill=black,label=below:$-3$] (k2) at (5,0) {};
			\node[draw,circle,fill=black,label=below:$-3$] (k3) at (5.5,0) {};
			
			\node[draw,circle,fill=black,label=below:$-3$] (l1) at (6,0) {};
			\node[draw,circle,fill=black,label=below:$-3$] (l2) at (6.5,0) {};
			\node[draw,circle,fill=black,label=below:$-3$] (l3) at (7,0) {};
			
			\node[draw,circle,fill=black,label=below:$-3$] (m1) at (7.5,0) {};
			\node[draw,circle,fill=black,label=below:$-3$] (m2) at (8,0) {};
			\node[draw,circle,fill=black,label=below:$-3$] (m3) at (8.5,0) {};
			
			\node[draw,circle,fill=black,label=below:$-3$] (n1) at (9,0) {};
			\node[draw,circle,fill=black,label=below:$-3$] (n2) at (9.5,0) {};
			\node[draw,circle,fill=black,label=below:$-3$] (n3) at (10,0) {};
			
			\node[draw,circle,fill=black,label=below:$-3$] (o1) at (10.5,0) {};
			\node[draw,circle,fill=black,label=below:$-3$] (o2) at (11,0) {};
			\node[draw,circle,fill=black,label=below:$-3$] (o3) at (11.5,0) {};

			  \path
				  (a) edge node[right]{}(b)
					(a) edge node[right]{}(c)
					
					(b) edge node[right]{}(d)
					(b) edge node[right]{}(e)
					
					(c) edge node[right]{}(f)
					(c) edge node[right]{}(g)
				
				  (d) edge node[right]{}(h)
					(d) edge node[right]{}(i)
					
					(e) edge node[right]{}(j)
					(e) edge node[right]{}(k)
					
					(f) edge node[right]{}(l)
					(f) edge node[right]{}(m)
					
					(g) edge node[right]{}(n)
					(g) edge node[right]{}(o)
				
				  (h) edge node[right]{}(h1)
					(h) edge node[right]{}(h2)
					(h) edge node[right]{}(h3)
					
          (i) edge node[right]{}(i1)
					(i) edge node[right]{}(i2)
					(i) edge node[right]{}(i3)
					
					(j) edge node[right]{}(j1)
					(j) edge node[right]{}(j2)
					(j) edge node[right]{}(j3)
					
					(k) edge node[right]{}(k1)
					(k) edge node[right]{}(k2)
					(k) edge node[right]{}(k3)
					
					(l) edge node[right]{}(l1)
					(l) edge node[right]{}(l2)
					(l) edge node[right]{}(l3)
					
					(m) edge node[right]{}(m1)
					(m) edge node[right]{}(m2)
					(m) edge node[right]{}(m3)
					
					(n) edge node[right]{}(n1)
					(n) edge node[right]{}(n2)
					(n) edge node[right]{}(n3)
					
					(o) edge node[right]{}(o1)
					(o) edge node[right]{}(o2)
					(o) edge node[right]{}(o3)
							
        ;
\end{tikzpicture}
\caption{$T_{G_{4}}(2,2,2,0| 0,0,0,3)$.}
       \label{cotree_example2}
\end{center}			
\end{figure}			

Notice that, we can directly apply Lemma \ref{lem2} at each node $\cup$ at level $3$ because its leaves have the same assignment $x_{4}=-3$.
And, after $m-1=3-1=2$ iterations in each node we obtain the cotree in Figure \ref{cotree_example3}.

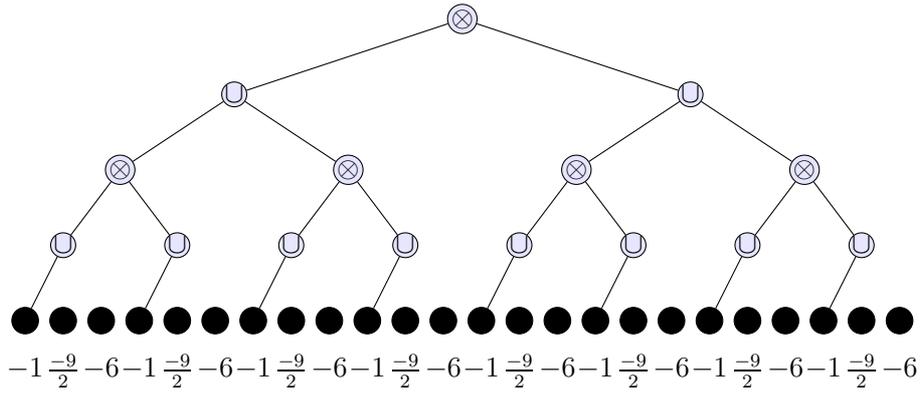
\begin{figure}[h!]
\begin{center}
\begin{tikzpicture}
  [scale=1,auto=left,every node/.style={circle,scale=0.9}]

                \node[draw, circle, fill=blue!10, inner sep=0] (a) at (5.75,4) {$\otimes$};
								
              \node[draw, circle, fill=blue!10, inner sep=0] (b) at (2.75,3) {$\cup$};
							 \node[draw, circle, fill=blue!10, inner sep=0] (c) at (8.75,3) {$\cup$};

        \node[draw, circle, fill=blue!10, inner sep=0] (d) at (1.25,2) {$\otimes$};
				\node[draw, circle, fill=blue!10, inner sep=0] (e) at (4.25,2) {$\otimes$};
				\node[draw, circle, fill=blue!10, inner sep=0] (f) at (7.25,2) {$\otimes$};
				\node[draw, circle, fill=blue!10, inner sep=0] (g) at (10.25,2) {$\otimes$};

       \node[draw, circle, fill=blue!10, inner sep=0] (h) at (0.5,1) {$\cup$};
       \node[draw, circle, fill=blue!10, inner sep=0] (i) at (2,1) {$\cup$};			
			\node[draw, circle, fill=blue!10, inner sep=0] (j) at (3.5,1) {$\cup$};
			\node[draw, circle, fill=blue!10, inner sep=0] (k) at (5,1) {$\cup$};
			\node[draw, circle, fill=blue!10, inner sep=0] (l) at (6.5,1) {$\cup$};
			\node[draw, circle, fill=blue!10, inner sep=0] (m) at (8,1) {$\cup$};
			\node[draw, circle, fill=blue!10, inner sep=0] (n) at (9.5,1) {$\cup$};
			\node[draw, circle, fill=blue!10, inner sep=0] (o) at (11,1) {$\cup$};
 			
			\node[draw,circle,fill=black,label=below:$-1$] (h1) at (0,0) {};
			\node[draw,circle,fill=black,label=below:$\frac{-9}{2}$] (h2) at (0.5,0) {};
			\node[draw,circle,fill=black,label=below:$-6$] (h3) at (1,0) {};
			
			\node[draw,circle,fill=black,label=below:$-1$] (i1) at (1.5,0) {};
			\node[draw,circle,fill=black,label=below:$\frac{-9}{2}$] (i2) at (2,0) {};
			\node[draw,circle,fill=black,label=below:$-6$] (i3) at (2.5,0) {};
			
			\node[draw,circle,fill=black,label=below:$-1$] (j1) at (3,0) {};
			\node[draw,circle,fill=black,label=below:$\frac{-9}{2}$] (j2) at (3.5,0) {};
			\node[draw,circle,fill=black,label=below:$-6$] (j3) at (4,0) {};
			
			\node[draw,circle,fill=black,label=below:$-1$] (k1) at (4.5,0) {};
			\node[draw,circle,fill=black,label=below:$\frac{-9}{2}$] (k2) at (5,0) {};
			\node[draw,circle,fill=black,label=below:$-6$] (k3) at (5.5,0) {};
			
			\node[draw,circle,fill=black,label=below:$-1$] (l1) at (6,0) {};
			\node[draw,circle,fill=black,label=below:$\frac{-9}{2}$] (l2) at (6.5,0) {};
			\node[draw,circle,fill=black,label=below:$-6$] (l3) at (7,0) {};
			
			\node[draw,circle,fill=black,label=below:$-1$] (m1) at (7.5,0) {};
			\node[draw,circle,fill=black,label=below:$\frac{-9}{2}$] (m2) at (8,0) {};
			\node[draw,circle,fill=black,label=below:$-6$] (m3) at (8.5,0) {};
			
			\node[draw,circle,fill=black,label=below:$-1$] (n1) at (9,0) {};
			\node[draw,circle,fill=black,label=below:$\frac{-9}{2}$] (n2) at (9.5,0) {};
			\node[draw,circle,fill=black,label=below:$-6$] (n3) at (10,0) {};
			
			\node[draw,circle,fill=black,label=below:$-1$] (o1) at (10.5,0) {};
			\node[draw,circle,fill=black,label=below:$\frac{-9}{2}$] (o2) at (11,0) {};
			\node[draw,circle,fill=black,label=below:$-6$] (o3) at (11.5,0) {};

			  \path
				  (a) edge node[right]{}(b)
					(a) edge node[right]{}(c)
					
					(b) edge node[right]{}(d)
					(b) edge node[right]{}(e)
					
					(c) edge node[right]{}(f)
					(c) edge node[right]{}(g)
				
				  (d) edge node[right]{}(h)
					(d) edge node[right]{}(i)
					
					(e) edge node[right]{}(j)
					(e) edge node[right]{}(k)
					
					(f) edge node[right]{}(l)
					(f) edge node[right]{}(m)
					
					(g) edge node[right]{}(n)
					(g) edge node[right]{}(o)
				
				  (h) edge node[right]{}(h1)
					
          (i) edge node[right]{}(i1)
					
					(j) edge node[right]{}(j1)
					
					(k) edge node[right]{}(k1)
					
					(l) edge node[right]{}(l1)
					
					(m) edge node[right]{}(m1)
					
					(n) edge node[right]{}(n1)
					
					(o) edge node[right]{}(o1)
							
        ;
\end{tikzpicture}
\caption{Iterations at level $4$.}
       \label{cotree_example3}
\end{center}			
\end{figure}			

Now, $16$ leaves are removed and have permanent values. The remaining $8$ leaves are relocated to level $3$ with value $x_{3}=-1$. Then, we move up to the cotree $T_{G_{3}}(2,2,0|0,0,2)$ as shown in Figure \ref{cotree_example4}.

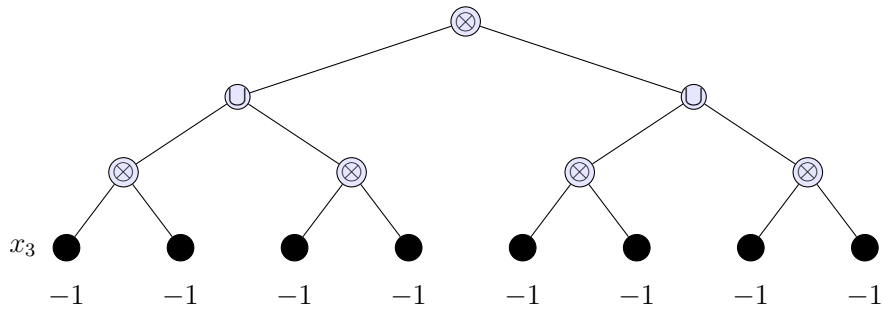
\begin{figure}[h!]
\begin{center}
\begin{tikzpicture}
  [scale=1,auto=left,every node/.style={circle,scale=0.9}]

                \node[draw, circle, fill=blue!10, inner sep=0] (a) at (5.75,4) {$\otimes$};
								
              \node[draw, circle, fill=blue!10, inner sep=0] (b) at (2.75,3) {$\cup$};
							 \node[draw, circle, fill=blue!10, inner sep=0] (c) at (8.75,3) {$\cup$};

        \node[draw, circle, fill=blue!10, inner sep=0] (d) at (1.25,2) {$\otimes$};
				\node[draw, circle, fill=blue!10, inner sep=0] (e) at (4.25,2) {$\otimes$};
				\node[draw, circle, fill=blue!10, inner sep=0] (f) at (7.25,2) {$\otimes$};
				\node[draw, circle, fill=blue!10, inner sep=0] (g) at (10.25,2) {$\otimes$};

       \node[draw,circle,fill=black,label=left:$x_{3}$,label=below:$-1$] (h) at (0.5,1) {};
       \node[draw,circle,fill=black,label=below:$-1$] (i) at (2,1) {};			
			\node[draw,circle,fill=black,label=below:$-1$] (j) at (3.5,1) {};
			\node[draw,circle,fill=black,label=below:$-1$] (k) at (5,1) {};
			\node[draw,circle,fill=black,label=below:$-1$] (l) at (6.5,1) {};
			\node[draw,circle,fill=black,label=below:$-1$] (m) at (8,1) {};
			\node[draw,circle,fill=black,label=below:$-1$] (n) at (9.5,1) {};
			\node[draw,circle,fill=black,label=below:$-1$] (o) at (11,1) {};
 									
			  \path
				  (a) edge node[right]{}(b)
					(a) edge node[right]{}(c)
					
					(b) edge node[right]{}(d)
					(b) edge node[right]{}(e)
					
					(c) edge node[right]{}(f)
					(c) edge node[right]{}(g)
				
				  (d) edge node[right]{}(h)
					(d) edge node[right]{}(i)
					
					(e) edge node[right]{}(j)
					(e) edge node[right]{}(k)
					
					(f) edge node[right]{}(l)
					(f) edge node[right]{}(m)
					
					(g) edge node[right]{}(n)
					(g) edge node[right]{}(o)
				
							
        ;
\end{tikzpicture}
\caption{$T_{G_{3}}(2,2,0|0,0,2)$.}
       \label{cotree_example4}
\end{center}			
\end{figure}

\newpage

Next, we can also directly apply Lemma \ref{lem1} at each node $\otimes$ at level $2$ because its leaves have the same assignment $x_{3}=-1$.
And, after $m-1=2-1=1$ iteration in each node we obtain the cotree in Figure \ref{cotree_example5}.

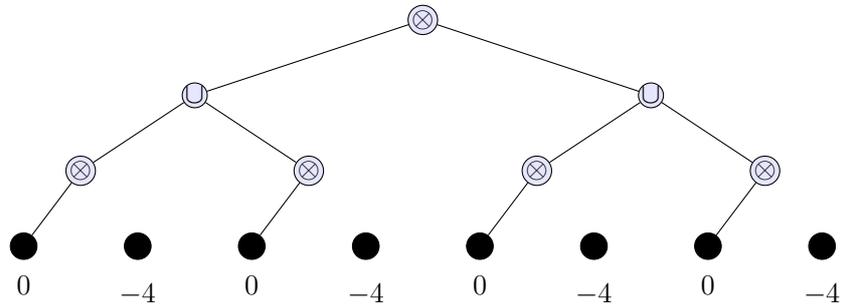
\begin{figure}[h!]
\begin{center}
\begin{tikzpicture}
  [scale=1,auto=left,every node/.style={circle,scale=0.9}]

                \node[draw, circle, fill=blue!10, inner sep=0] (a) at (5.75,4) {$\otimes$};
								
              \node[draw, circle, fill=blue!10, inner sep=0] (b) at (2.75,3) {$\cup$};
							 \node[draw, circle, fill=blue!10, inner sep=0] (c) at (8.75,3) {$\cup$};

        \node[draw, circle, fill=blue!10, inner sep=0] (d) at (1.25,2) {$\otimes$};
				\node[draw, circle, fill=blue!10, inner sep=0] (e) at (4.25,2) {$\otimes$};
				\node[draw, circle, fill=blue!10, inner sep=0] (f) at (7.25,2) {$\otimes$};
				\node[draw, circle, fill=blue!10, inner sep=0] (g) at (10.25,2) {$\otimes$};

       \node[draw,circle,fill=black,label=below:$0$] (h) at (0.5,1) {};
       \node[draw,circle,fill=black,label=below:$-4$] (i) at (2,1) {};			
			\node[draw,circle,fill=black,label=below:$0$] (j) at (3.5,1) {};
			\node[draw,circle,fill=black,label=below:$-4$] (k) at (5,1) {};
			\node[draw,circle,fill=black,label=below:$0$] (l) at (6.5,1) {};
			\node[draw,circle,fill=black,label=below:$-4$] (m) at (8,1) {};
			\node[draw,circle,fill=black,label=below:$0$] (n) at (9.5,1) {};
			\node[draw,circle,fill=black,label=below:$-4$] (o) at (11,1) {};
 									
			  \path
				  (a) edge node[right]{}(b)
					(a) edge node[right]{}(c)
					
					(b) edge node[right]{}(d)
					(b) edge node[right]{}(e)
					
					(c) edge node[right]{}(f)
					(c) edge node[right]{}(g)
				
				  (d) edge node[right]{}(h)

					(e) edge node[right]{}(j)
					
					(f) edge node[right]{}(l)
					
					(g) edge node[right]{}(n)
									
							
        ;
\end{tikzpicture}
\caption{Iterations at level $3$.}
       \label{cotree_example5}
\end{center}			
\end{figure}

Now, $4$ leaves are removed and have permanent values. The remaining $4$ leaves are relocated to level $2$ with value $x_{2}=0$. Then, we move up to the cotree $T_{G_{2}}(2,0|0,2)$ as shown in Figure \ref{cotree_example6}.

\begin{figure}[h!]
\begin{center}
\begin{tikzpicture}
  [scale=1,auto=left,every node/.style={circle,scale=0.9}]

                \node[draw, circle, fill=blue!10, inner sep=0] (a) at (5.75,4) {$\otimes$};
								
              \node[draw, circle, fill=blue!10, inner sep=0] (b) at (2.75,3) {$\cup$};
							 \node[draw, circle, fill=blue!10, inner sep=0] (c) at (8.75,3) {$\cup$};

        \node[draw,circle,fill=black,label=left:$x_{2}$,label=below:$0$] (d) at (1.25,2) {};
				\node[draw,circle,fill=black,label=below:$0$] (e) at (4.25,2) {};
				\node[draw,circle,fill=black,label=below:$0$] (f) at (7.25,2) {};
				\node[draw,circle,fill=black,label=below:$0$] (g) at (10.25,2) {};

			  \path
				  (a) edge node[right]{}(b)
					(a) edge node[right]{}(c)
					
					(b) edge node[right]{}(d)
					(b) edge node[right]{}(e)
					
					(c) edge node[right]{}(f)
					(c) edge node[right]{}(g)

							
        ;
\end{tikzpicture}
\caption{$T_{G_{2}}(2,0|0,2)$.}
       \label{cotree_example6}
\end{center}			
\end{figure}
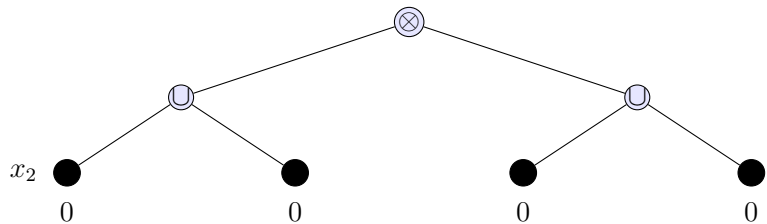
Hence, we apply \subiib~ in each pair of leaves at each node $\cup$ at level $1$. And, after $1$ iteration in each node we obtain the cotree in Figure \ref{cotree_example7}

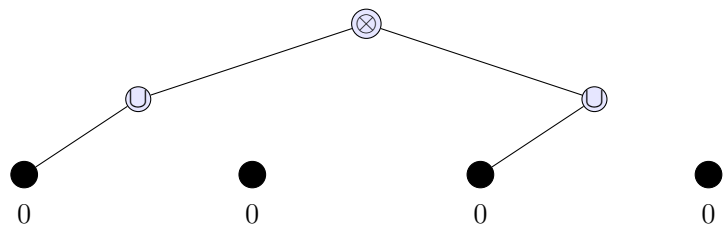
\begin{figure}[h!]
\begin{center}
\begin{tikzpicture}
  [scale=1,auto=left,every node/.style={circle,scale=0.9}]

                \node[draw, circle, fill=blue!10, inner sep=0] (a) at (5.75,4) {$\otimes$};
								
              \node[draw, circle, fill=blue!10, inner sep=0] (b) at (2.75,3) {$\cup$};
							 \node[draw, circle, fill=blue!10, inner sep=0] (c) at (8.75,3) {$\cup$};

        \node[draw,circle,fill=black,label=below:$0$] (d) at (1.25,2) {};
				\node[draw,circle,fill=black,label=below:$0$] (e) at (4.25,2) {};
				\node[draw,circle,fill=black,label=below:$0$] (f) at (7.25,2) {};
				\node[draw,circle,fill=black,label=below:$0$] (g) at (10.25,2) {};

			  \path
				  (a) edge node[right]{}(b)
					(a) edge node[right]{}(c)
					
					(b) edge node[right]{}(d)
									
					(c) edge node[right]{}(f)
								  								
							
        ;
\end{tikzpicture}
\caption{Iterations at level $2$.}
       \label{cotree_example7}
\end{center}			
\end{figure}

Now, $2$ leaves are removed and have permanent values. The remaining $2$ leaves are relocated to level $1$ with value $x_{1}=0$. Then, we move up to the cotree $T_{G_{1}}(0|2)$ as shown in Figure \ref{cotree_example8}.

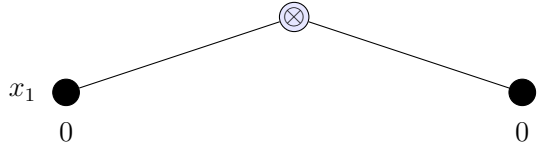
\begin{figure}[h!]
\begin{center}
\begin{tikzpicture}
  [scale=1,auto=left,every node/.style={circle,scale=0.9}]

                \node[draw, circle, fill=blue!10, inner sep=0] (a) at (5.75,4) {$\otimes$};
								
              \node[draw,circle,fill=black,label=left:$x_{1}$,label=below:$0$] (b) at (2.75,3) {};
			  \node[draw,circle,fill=black,label=below:$0$] (c) at (8.75,3) {};

			  \path
				  (a) edge node[right]{}(b)
				  (a) edge node[right]{}(c)
					
        ;
\end{tikzpicture}
\caption{$T_{G_{1}}(0|2)$.}
       \label{cotree_example8}
\end{center}			
\end{figure}

Finally, we apply Lemma \ref{lem1} at the root $\otimes$. And, after $1$ iteration we obtain the cotree in Figure \ref{cotree_example9}.

\begin{figure}[h!]
\begin{center}
\begin{tikzpicture}
  [scale=1,auto=left,every node/.style={circle,scale=0.9}]

                \node[draw, circle, fill=blue!10, inner sep=0] (a) at (5.75,4) {$\otimes$};
								
              \node[draw,circle,fill=black,label=below:$\frac{1}{2}$] (b) at (2.75,3) {};
			  \node[draw,circle,fill=black,label=below:$-2$] (c) at (8.75,3) {};

			  \path
				  (a) edge node[right]{}(b)
				  				
        ;
\end{tikzpicture}
\caption{Iterations at level $1$.}
       \label{cotree_example9}
\end{center}
\end{figure}
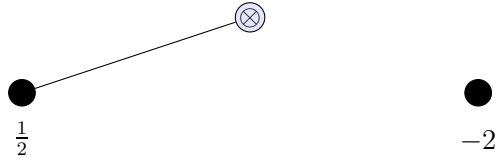

Therefore, by Theorem \ref{main1}, $3$ is an eigenvalue of the cograph with multiplicity $2$. $21$ eigenvalues are less than $3$ and $1$ eigenvalue is greater than $3$. And, using Theorem \ref{main2} we know $0$ is an eigenvalue of multiplicity $16$.

\end{Ex}

\begin{Remark}
In our notation, we apply Algorithm Diagonalize to a balanced cotree $T_{G_{r}}$ $(a_{1},\ldots,a_{r-1},0|0,\ldots,0,a_{r})$ with value $x_{r}$ in the leaves. Then, we proceed with the iterations at level $r$. Once the level is done. We have removed $a_{1}a_{2}\cdots a_{r-1}(a_{r}-1)$ leaves and the remaining leaves are relocated to level $r-1$ with value $x_{r-1}$. And, we move up to the cotree $T_{G_{r-1}}$ $(a_{1},\ldots,a_{r-2},0|0,\ldots,0,a_{r-1})$ and so on.
\end{Remark}

In the next section, we use an inverse procedure to construct integral cographs. It means that, we start at some level $i$ with the cotree $T_{G_{i}}$ $(a_{1},\ldots,a_{i-1},0|0,\ldots,0,a_{i})$ and we proceed from top to bottom until we compute the input $x$ in the cotree $T_{G_{r}}$ $(a_{1},\ldots,a_{r-1},0|0,\ldots,0,a_{r})$. Then, we will have the eigenvalue $-x$.


\section{Recurrence Relation}
\label{recurrence}

In this section, we start at some level $i$ in Algorithm Diagonalize with the cotree $T_{G_{i} }$ and the value $x_{i}$ in the leaves. Then, we proceed to the previous steps, it means that, we go from the cotree $T_{G_{i}}(a_{1},\ldots,a_{i-1},0|$ $0,\ldots,0,a_{i})$ to $T_{G_{i+1}}$ $(a_{1},\ldots,a_{i},0|0,\ldots,0,a_{i+1})$. From $T_{G_{i+1}}(a_{1},\ldots,a_{i},0|0,\ldots,0,a_{i+1})$ to $T_{G_{i+2}}(a_{1},\ldots,a_{i+2},0|0,\ldots,0,a_{i+1})$ and so on. Until we compute the inputs $x$ and
$T_{G_{r}}(a_{1},\ldots,a_{r-1},0|0,\ldots,0,a_{r})$. We show this process gives a recurrence relation, whose solution leads to an integer eigenvalue of the cograph $G$.

\begin{Lem}
\label{01}
Let $G$ be a  cograph with balanced cotree $T_{G_{r}}$ $(a_1,\ldots,a_{r-1},0|$ $0,\ldots,0,a_r)$ of order $n =a_1 a_2 \ldots a_{r-1} a_r$.
The Algorithm Diagonalize with input $(T_{G_{r}}, x)$ produces leaves $x_{i}$ in the cotree $T_{G_{i}}  (a_1,  \ldots,a_{i-1}, 0| 0, \ldots,0, a_i)$ with value zero at level $i$ if $i$ is even, and with value $1$ if $i$ is odd, for $ 1 \leq i \leq r-1$  if and only if the balanced cotree $T_{G_{i+1}}(a_1,  \ldots,a_{i}, 0| 0, \ldots,0, a_{i+1})$, in the previous step,
 has  initialized with value
$$
x_{i+1}= \left\{ \begin{array}{rl}
 1-a_{i+1}  &\mbox{ if $i$ is even} \\
 a_{i+1} &\mbox{ if $i$ is odd}
       \end{array} \right.
$$
for $ 1 \leq i \leq r-1$. And, at the iterations at level $i+1$, the removed leaves have non-null permanent values.
\end{Lem}
{\bf Proof:} Given $1\leq i\leq r-1$. First, we consider $i$ even and the cotree $T_{G_{i}}  (a_1,  \ldots,a_{i-1}, 0| 0, \ldots,0,a_{i})$ with each leave having value $x_{i}=0$ as in Figure \ref{fig00}.

Now, we consider we are in the previous level with the cotree $T_{G_{i+1}}$ $(a_1,  \ldots,a_{i}, 0| 0, \ldots,0, a_{i+1})$ having leaves with value $x_{i+1}$.

 If $x_{i+1}=1$, then at each node $\otimes$ at level $i$ we can apply \subib~ at its leaves $m-1=a_{i+1}-1$ times. At the last iteration, at each node, the remaining leave has the value $$x_{i}=d_{l}\leftarrow 1.$$
Contradiction, because $x_{i}=0$. Therefore, $x_{i+1}\neq 1$.

 Since the leaves have value $x_{i+1}\neq 1$, then we can apply Lemma \ref{lem1}. Hence, we perform $m-1=a_{i+1}-1$ iterations at each node $\otimes$ at level $i$. At the last iteration, at each node, the remaining leave receives the following assignment
$$x_{i}=d_{l}\leftarrow\frac{x_{i+1}+a_{i+1}-1}{a_{i+1}-1+1}=\frac{x_{i+1}+a_{i+1}-1}{a_{i+1}}.$$
But, $x_{i}=0$ if and only if $x_{i+1}+a_{i+1}-1=0$. It implies $x_{i+1}=1-a_{i+1}$.
And, the removed leaves at each node $\otimes$ at level $i$ receive
$$d_{k}\leftarrow\frac{j+1}{j}(x_{i+1}-1)\neq 0,$$
for $1\leq j\leq a_{i+1}-1$.

\begin{figure}[h!]
\begin{center}
\begin{tikzpicture}
  [scale=1,auto=left,every node/.style={circle,scale=0.9}]

       \node[draw, circle, fill=blue!10, inner sep=0,label=left:$a_{i-1}$,label=above:$T_{G_{i}}$] (a) at (1,1) {$\cup$};
 			 \node[draw,circle,fill=black, label=below:$x_{i}$,label=left:$a_{i}$] (b) at (0,0) {};
			\node[draw,circle,fill=black, label=below:$x_{i}$] (c) at (2,0) {};

				\node[draw, circle, fill=blue!10, inner sep=0,label=left:$a_{i}$,label=above:$T_{G_{i+1}}$] (d) at (6,1) {$\otimes$};
 			 \node[draw,circle,fill=black, label=below:$x_{i+1}$,label=left:$a_{i+1}$] (e) at (5,0) {};
			\node[draw,circle,fill=black, label=below:$x_{i+1}$] (f) at (7,0) {};     			
							
							\draw[dashed, very thick] (0.3,0) -- (1.7,0);
							\draw[dashed, very thick] (5.3,0) -- (6.7,0);
							\draw [->] (3,0.5) -- (4,0.5);
  \path
          (a) edge node[right]{}(b)
          (a) edge node[right]{}(c)
							
					(d) edge node[right]{}(e)
          (d) edge node[right]{}(f)
        ;
\end{tikzpicture}
\caption{$i$ even.}
       \label{fig00}
			
\end{center}			
\end{figure}

Now, $i$ is odd and the cotree $T_{G_{i}} (a_1,  \ldots,a_{i-1}, 0| 0, \ldots,0,a_{i})$ has leaves with value $x_{i}=1$ as in Figure \ref{fig01}.

If $x_{i+1}=0$, then at each node $\cup$ at level $i$ we can apply \subiib~ at its leaves $m-1=a_{i+1}-1$ times. At the last iteration, at each node, the remaining leave has value

$$x_{i}=d_{l}\leftarrow 0.$$

Contradiction, because $x_{i}=1$. Therefore, $x_{i+1}\neq 0$.

Since $x_{i+1}\neq 0$, then we can apply Lemma \ref{lem2} at each node $\cup$. We perform $m-1=a_{i+1}-1$ iterations at each node and at the last iteration, in each node, the remaining leave receives

$$x_{i}=d_{l}\leftarrow\frac{x_{i+1}}{a_{i+1}-1+1}.$$

It implies $x_{i+1}=a_{i+1}\cdot x_{i}=a_{i+1}$.
And, the removed leaves at each node $\cup$ have permanent values
$$d_{k}\leftarrow\frac{j+1}{j}(x_{i+1})\neq 0,$$
for $1\leq j\leq a_{i+1}-1$.

\begin{figure}[h!]
\begin{center}
\begin{tikzpicture}
  [scale=1,auto=left,every node/.style={circle,scale=0.9}]

       \node[draw, circle, fill=blue!10, inner sep=0,label=left:$a_{i-1}$,label=above:$T_{G_{i}}$] (a) at (1,1) {$\otimes$};
 			 \node[draw,circle,fill=black, label=below:$x_{i}$,label=left:$a_{i}$] (b) at (0,0) {};
			\node[draw,circle,fill=black, label=below:$x_{i}$] (c) at (2,0) {};

				\node[draw, circle, fill=blue!10, inner sep=0,label=left:$a_{i}$,label=above:$T_{G_{i+1}}$] (d) at (6,1) {$\cup$};
 			 \node[draw,circle,fill=black, label=below:$x_{i+1}$,label=left:$a_{i+1}$] (e) at (5,0) {};
			\node[draw,circle,fill=black, label=below:$x_{i+1}$] (f) at (7,0) {};     			
							
							\draw[dashed, very thick] (0.3,0) -- (1.7,0);
							\draw[dashed, very thick] (5.3,0) -- (6.7,0);
							\draw [->] (3,0.5) -- (4,0.5);
  \path
          (a) edge node[right]{}(b)
          (a) edge node[right]{}(c)
							
					(d) edge node[right]{}(e)
          (d) edge node[right]{}(f)
        ;
\end{tikzpicture}
\caption{$i$ odd.}
       \label{fig01}
\end{center}			
\end{figure}			

Given $1\leq i\leq r-1$. If $i$ is even and the cotree $T_{G_{i+1}}  (a_1,  \ldots,a_{i}, 0| 0, \ldots,0, a_{i+1})$ has leaves with value $x_{i+1}=1-a_{i+1}$. As $a_{i+1}\geq 2$, we have $x_{i+1}\neq 1$. Then, we can use Lemma \ref{lem1} at each node $\otimes$, at level $i$, $(a_{i+1}-1)$ times. After the last iteration at each node the remaining leave receives
$$x_{i}=d_{l}\leftarrow\frac{x_{i+1}+a_{i+1}-1}{a_{i+1}-1+1}=\frac{1-a_{i+1}+a_{i+1}-1}{a_{i+1}}=0.$$
And, the removed leaves at each node receive
$$d_{k}\leftarrow\frac{j+1}{j}(x_{i+1}-1)=\frac{j+1}{j}(1-a_{i+1}-1)=\frac{j+1}{j}(-a_{i+1})\neq 0,$$
 for $1\leq j\leq a_{i+1}-1$.

And, if $i$ is odd. Then $x_{i+1}=a_{i+1}\geq 2$ and we can use Lemma \ref{lem2} at each node $\cup$ at level $i$, $(a_{i+1}-1)$ times. After the last iteration, each remaining leave at level $i+1$ receives
$$x_{i}=d_{l}\leftarrow\frac{a_{i+1}}{a_{i+1}-1+1}=1.$$
And, the removed leaves at each node receive
$$d_{k}\leftarrow\frac{j+1}{j}(x_{i+1})=\frac{j+1}{j}(a_{i+1})\neq 0,$$
 for $1\leq j\leq a_{i+1}-1$.

The next Lemma guarantees we can use Lemmas \ref{lem1} and \ref{lem2} recursively, in the opposite way.

\begin{Lem}
\label{lema00}
Let $G$ be a  cograph with balanced cotree $T_{G_{r}}$ $(a_1,\ldots,a_{r-1}, 0$ $|0,\ldots,0,a_r)$ of order $n =a_1 a_2 \ldots a_{r-1} a_r$.
If the Algorithm Diagonalize with input $(T_{G_{r}}, x)$ produces leaves $x_{i}$ at level $i$ with value zero if $i$ is even, or with value $1$ if $i$ is odd, for $ 1 \leq i \leq r-1$.
Then, in the previous steps, the algorithm has produced leaves in the cotree $T_{G_{i+j}}$ with value $x_{i+j}$ for $j\geq 1$ such that

$$
x_{i+j} \left\{ \begin{array}{rl}
 \neq 0, & \\
 \neq 1, & \\
 \in\mathbb{Z}. & \\
       \end{array} \right.
$$

And, at each level, the removed leaves have non-null permanent.
\end{Lem}
{\bf Proof:} We prove by induction on $j$. We start with $j=1$. If $i$ is even then $x_{i}=0$ and by Lemma \ref{01} we have $x_{i+1}=1-a_{i+1}$. As $a_{i+1}\geq 2$ then $x_{i+1}\leq -1$. And $x_{i+1}=1-a_{i+1}\in\mathbb{Z}$ because $a_{i+1}\in\mathbb{Z}$.

In this case, the removed leaves receive $$d_{k}\leftarrow\frac{j+1}{j}(1-a_{i+1}-1) =\frac{j+1}{j}(-a_{i+1})\neq 0,$$
for $1\leq j\leq a_{i+1}-1$.

If $i$ is odd then $x_{i}=1$ and by Lemma \ref{01} we have $x_{i+1}=a_{i+1}$. As $a_{i+1}\in\mathbb{Z}$ and $a_{i+1}\geq 2$.
In this case, the removed leaves receive $$d_{k}\leftarrow\frac{j+1}{j}(x_{i+1})=\frac{j+1}{j}(a_{i+1})\neq 0,$$
for $1\leq j\leq a_{i+1}-1$.

Now, suppose $x_{i+j}\neq 0$, $\neq 1$ and it is an integer number.

First, we consider $i+j$ even as shown in Figure \ref{fig02}.
\begin{figure}[h!]
\begin{center}
\begin{tikzpicture}
  [scale=1,auto=left,every node/.style={circle,scale=0.9}]

       \node[draw, circle, fill=blue!10, inner sep=0,label=left:$a_{i+j-1}$,label=above:$T_{G_{i+j}}$] (a) at (1,1) {$\cup$};
 			 \node[draw,circle,fill=black, label=below:$x_{i+j}$,label=left:$a_{i+j}$] (b) at (0,0) {};
			\node[draw,circle,fill=black, label=below:$x_{i+j}$] (c) at (2,0) {};

				\node[draw, circle, fill=blue!10, inner sep=0,label=left:$a_{i+j}$,label=above:$T_{G_{i+j+1}}$] (d) at (6,1) {$\otimes$};
 			 \node[draw,circle,fill=black, label=below:$x_{i+j+1}$,label=left:$a_{i+j+1}$] (e) at (5,0) {};
			\node[draw,circle,fill=black, label=below:$x_{i+j+1}$] (f) at (7,0) {};     			
							
							\draw[dashed, very thick] (0.3,0) -- (1.7,0);
							\draw[dashed, very thick] (5.3,0) -- (6.7,0);
							\draw [->] (3,0.5) -- (4,0.5);
  \path
          (a) edge node[right]{}(b)
          (a) edge node[right]{}(c)
							
					(d) edge node[right]{}(e)
          (d) edge node[right]{}(f)
        ;
\end{tikzpicture}
\caption{$i$ even.}
       \label{fig02}
\end{center}
\end{figure}

If $x_{i+j+1}=0$ then we can use Lemma \ref{lem1}. And, at the last iteration $j=a_{i+j+1}-1$, we have
$$x_{i+j}=d_{l}\leftarrow\frac{x_{i+j+1}+a_{i+j+1}-1}{a_{i+j+1}-1+1}=\frac{a_{i+j+1}-1}{a_{i+j+1}}\notin\mathbb{Z}.$$ As it is a contradiction we have $x_{i+j+1}\neq 0$.

If $x_{i+j+1}=1$ then we apply \subib~ $a_{i+j+1}-1$ times. And, at the last iteration, we have
$$x_{i+j}=d_{l}\leftarrow 1.$$
And, as it contradicts the hypothesis, $x_{i+j+1}\neq 1$.

We already know that $x_{i+j+1}\neq 1$. Then, we can use Lemma \ref{lem1} to compute $x_{i+j}$. Notice that, at the last iteration at the node we have $j=a_{i+j+1}-1$ and it implies
$$x_{i+j}=d_{l}\leftarrow\frac{x_{i+j+1}+a_{i+j+1}-1}{a_{i+j+1}-1+1},$$
hence $x_{i+j+1}=a_{i+j+1}\cdot x_{i+j}-a_{i+j+1}+1\in\mathbb{Z}$.

 And, the removed leaves receive $$d_{k}\leftarrow\frac{j+1}{j}(x_{i+j+1}-1)=\frac{j+1}{j}(a_{i+j+1})(x_{i+j}-1)\neq 0,$$
 for $1\leq j\leq a_{i+j+1}-1$.

Now, we consider $i+j$ odd as shown in Figure \ref{fig03}.
\begin{figure}[h!]
\begin{center}
\begin{tikzpicture}
  [scale=1,auto=left,every node/.style={circle,scale=0.9}]

       \node[draw, circle, fill=blue!10, inner sep=0,label=left:$a_{i+j-1}$,label=above:$T_{G_{i+j}}$] (a) at (1,1) {$\otimes$};
 			 \node[draw,circle,fill=black, label=below:$x_{i+j}$,label=left:$a_{i+j}$] (b) at (0,0) {};
			\node[draw,circle,fill=black, label=below:$x_{i+j}$] (c) at (2,0) {};

				\node[draw, circle, fill=blue!10, inner sep=0,label=left:$a_{i+j}$,label=above:$T_{G_{i+j+1}}$] (d) at (6,1) {$\cup$};
 			 \node[draw,circle,fill=black, label=below:$x_{i+j+1}$,label=left:$a_{i+j+1}$] (e) at (5,0) {};
			\node[draw,circle,fill=black, label=below:$x_{i+j+1}$] (f) at (7,0) {};     			
							
							\draw[dashed, very thick] (0.3,0) -- (1.7,0);
							\draw[dashed, very thick] (5.3,0) -- (6.7,0);
							\draw [->] (3,0.5) -- (4,0.5);
  \path
          (a) edge node[right]{}(b)
          (a) edge node[right]{}(c)
							
					(d) edge node[right]{}(e)
          (d) edge node[right]{}(f)
        ;
\end{tikzpicture}
\caption{$i$ odd.}
       \label{fig03}
\end{center}
\end{figure}	
		
If $x_{i+j+1}=0$ then we can use \subiib~ $(a_{i+j+1}-1)$ times at the node $\cup$. And, at the last iteration, we  have
$$x_{i+j}=d_{l}\leftarrow 0.$$
As it is a contradiction, we have $x_{i+j+1}\neq 0$.

Now, if $x_{i+j+1}=1$ then we can use Lemma \ref{lem2}. And, at the last iteration $j=a_{i+j+1}-1$, we obtain $$x_{i+j}=d_{l}\leftarrow\frac{x_{i+j+1}}{a_{i+j+1}-1+1}=\frac{1}{a_{i+j+1}}\notin\mathbb{Z}.$$
Therefore, $x_{i+j+1}\neq 1$.

Now, we already know that $x_{i+j+1}\neq 0$. Then, we can use Lemma \ref{lem2}. And, at the last iteration $j=a_{i+j+1}-1$, we have
$$x_{i+j}=d_{l}\leftarrow\frac{x_{i+j+1}}{a_{i+j+1}-1+1}.$$
It implies that $x_{i+j+1}=a_{i+j+1}\cdot x_{i+j}\in\mathbb{Z}$ because $x_{i+j}$, $a_{i+j+1}\in\mathbb{Z}$.

In this case, the removed leaves receive $$d_{k}\leftarrow\frac{j+1}{j}(x_{i+j+1})=\frac{j+1}{j}(a_{i+j+1}\cdot x_{i+j})\neq 0,$$
for $1\leq j\leq a_{i+j+1}-1$.


Suppose we have the cotree $T_{G_{i}}(a_{1},\ldots,a_{i-1},0|0,\ldots,0,a_{i})$ for some $i\in\{1,\ldots,r-1\}$ and its leaves have value $x_{i}=1$, if $i$ is odd, and   $x_{i}=0$ if $i$ is even.
Then, Lemma \ref{lema00} guarantees we can recursively use Lemmas \ref{lem1} and \ref{lem2} to retrieve the value of the input $x$ in Algorithm Diagonalize. It means that, we can compute eigenvalues of the cograph $G$ that has cotree $T_{G_{r}}(a_{1},\ldots,a_{r-1},0|0,\ldots,0,a_{r})$ using the following recurrence relation.

\[ \begin{cases}
      x_{i}=1 & \mbox{if $i$ is odd,} \\
		  x_{i}=0 & \mbox{if $i$ is even,} \\
      x_{i+j+1}=a_{i+j+1}\cdot x_{i+j} & \mbox{if $i+j$ is odd,} \\
			x_{i+j+1}=a_{i+j+1}\cdot x_{i+j}-a_{i+j+1}+1 & \mbox{if $i+j$ is even.}
   \end{cases}
\]

Now, we solve the recurrence relation above. For this, we define the following parameter.

\begin{Def}
Let $a = (a_1, a_2, \ldots, a_{n}) $  be a fixed sequence of positive integers. We define the following parameter

  $$   \gamma _{n,l}=  \left\{\begin{array}{ccc}
            a_n a_{n-1} a_{n-2} \ldots a_{l}& if &  1 \leq l \leq n-1 \\
                a_n  & if &    l=n .     \\
   \end{array}\right.$$
 \end{Def}

In the next two theorems we start with the cotree $T_{G_{i}}$ $(a_{1},\ldots,a_{i-1},0|$ $0,\ldots,0,a_{i})$ and we proceed from top to button.

\begin{Thr}
\label{rec1}
The recurrence relation for $i$ even is
\[ \begin{cases}
      x_{i}=0 & \mbox{if $i$ is even,} \\
      x_{i+j}=a_{i+j}\cdot x_{i+j-1} & \mbox{if $i+j-1$ is odd,} \\
			x_{i+j}=a_{i+j}\cdot x_{i+j-1}-a_{i+j}+1 & \mbox{if $i+j-1$ is even.}
   \end{cases}
\]
The solution is
\[
 \begin{cases}
      x_{i+j}=\sum_{k=1}^{j}\gamma_{i+j,i+k}(-1)^{k}+1 & \mbox{if $i+j$ is odd,} \\
      x_{i+j}=\sum_{k=1}^{j}\gamma_{i+j,i+k}(-1)^{k} & \mbox{if $i+j$ is even.}
    \end{cases}
\]
\end{Thr}
{\bf Proof:} We prove by induction on $j\geq 1$. For the basis of induction we consider $j=1$.
Using the recurrence relation for $i$ even we obtain
$$x_{i+1}=a_{i+1}\cdot x_{i}-a_{i+1}+1=-a_{i+1}+1.$$
And, $i+j=i+1$ is odd. Then, using the recurrence formula we can compare with
$$x_{i+1}=\sum_{k=1}^{1}\gamma_{i+1,i+k}(-1)^{k}+1=\gamma_{i+1,i+1}(-1)^{1}+1=-a_{i+1}+1.$$
The basis is done.

Now we have two cases: $i+j$ is even or odd.

We start with $i+j$ even. In this case $i+j-1$ is odd and by the recurrence relation we have
\begin{eqnarray}
\nonumber 
  x_{i+j} &=& a_{i+j}\cdot x_{i+j-1}= a_{i+j}(\sum_{k=1}^{j-1}\gamma_{i+j-1,i+k}(-1)^{k}+1)\\
\nonumber
          &=& \sum_{k=1}^{j-1}\gamma_{i+j,i+k}(-1)^{k}+a_{i+j}=
\sum_{k=1}^{j-1}\gamma_{i+j,i+k}(-1)^{k}+\gamma_{i+j,i+j}(-1)^{j},
\end{eqnarray}

and $(-1)^{j}=1$ because $j$ is even. It implies that

$$x_{i+j}=\sum_{k=1}^{j}\gamma_{i+j,i+k}(-1)^{k}.$$
Now, we consider $i+j$ odd. In this case $i+j-1$ is even and by the recurrence relation we have

\begin{eqnarray}
 \nonumber 
   x_{i+j}&=& a_{i+j}\cdot x_{i+j-1}-a_{i+j}+1=a_{i+j}(\sum_{k=1}^{j-1}\gamma_{i+j-1,i+k}(-1)^{k})-a_{i+j}+1 \\
  \nonumber
   &=& \sum_{k=1}^{j-1}\gamma_{i+j,i+k}(-1)^{k}-a_{i+j}+1,
\end{eqnarray}

and $-a_{i+j}=-\gamma_{i+j,i+j}=\sum_{k=j}^{j}\gamma_{i+j,i+k}(-1)^{j}$. As $(-1)^{j}=-1$ because $j$ is odd, we have
$$x_{i+j}=\sum_{k=1}^{j}\gamma_{i+j,i+k}(-1)^{k}+1.$$

\begin{Thr}
\label{rec2}
The recurrence relation for $i$ odd is
\[ \begin{cases}
      x_{i}=0 & \mbox{if $i$ is even,} \\
      x_{i+j}=a_{i+j}\cdot x_{i+j-1} & \mbox{if $i+j-1$ is odd,} \\
			x_{i+j}=a_{i+j}\cdot x_{i+j-1}-a_{i+j}+1 & \mbox{if $i+j-1$ is even.}
   \end{cases}
\]
The solution is
\[ \begin{cases}
      x_{i+j}=\sum_{k=1}^{j}\gamma_{i+j,i+k}(-1)^{k+1}+1 & \mbox{if $i+j$ is odd,} \\
		  x_{i+j}=\sum_{k=1}^{j}\gamma_{i+j,i+k}(-1)^{k+1} & \mbox{if $i+j$ is even.}
    \end{cases}
\]
\end{Thr}
{\bf Proof:} We prove by induction on $j\geq 1$. For the basis of induction we consider $j=1$.
Using the recurrence relation for $i$ odd we obtain
$$x_{i+1}=a_{i+1}\cdot x_{i}=a_{i+1}.$$
And, $i+j=i+1$ is even. Then using the recurrence formula we can compare with
$$x_{i+1}=\sum_{k=1}^{1}\gamma_{i+1,i+k}(-1)^{k+1}=\gamma_{i+1,i+1}(-1)^{2}=a_{i+1}.$$
The basis is done.

Now we have two cases: $i+j$ is even or odd.

We start with $i+j$ even. In this case $i+j-1$ is odd and by the recurrence relation we have

\begin{eqnarray}
 \nonumber 
  x_{i+j} &=& a_{i+j}\cdot x_{i+j-1}=a_{i+j}(\sum_{k=1}^{j-1}\gamma_{i+j-1,i+k}(-1)^{k+1}+1) \\
\nonumber
          &=& \sum_{k=1}^{j-1}\gamma_{i+j,i+k}(-1)^{k+1}+a_{i+j}.
\end{eqnarray}

And $(-1)^{j+1}=1$ because $j$ is odd. Then $a_{i+j}=(-1)^{j+1}\gamma_{i+j,i+j}$. It implies that

$$x_{i+j}=\sum_{k=1}^{j}\gamma_{i+j,i+k}(-1)^{k+1}.$$

Now, we consider $i+j$ odd. In this case $i+j-1$ is even and by the recurrence relation we have

\begin{eqnarray}
\nonumber 
  x_{i+j} &=& a_{i+j}\cdot x_{i+j-1}-a_{i+j}+1=a_{i+j}(\sum_{k=1}^{j-1}\gamma_{i+j-1,i+k}(-1)^{k+1})-a_{i+j}+1 \\
\nonumber
   &=& \sum_{k=1}^{j-1}\gamma_{i+j,i+k}(-1)^{k+1}+(-1)^{j+1}\gamma_{i+j,i+j}+1,
\end{eqnarray}

because $j$ is even. Then,

$$x_{i+j}=\sum_{k=1}^{j}\gamma_{i+j,i+k}(-1)^{k+1}+1.$$

Using Theorems \ref{rec1} and \ref{rec2} we can retrieve the eigenvalues of $T_{G_{r}}$ $(a_{1},\ldots,$ $a_{r-1},0|$ $0,\ldots,0,a_{r})$ that generate $x_{i}=1$, if $i$ is odd, and $x_{i}=0$, if $i$ is even, in  the cotree $T_{G_{i}}(a_{1},\ldots,a_{i-1},0|0,\ldots,0,a_{i})$.

From now on, we denote the input $x_{r}$ of the Algorithm Diagonalize to the cotree  $T_{G_{r}}(a_{1},\ldots,a_{r-1},0|0,\ldots,0,a_{r})$ that we retrieve from level $i$ by $x_{r}^{i}$.
The next two corollaries follow directly from Theorems \ref{rec1} and \ref{rec2} when $i+j=r$.
\begin{Cor}
\label{correc1}
If $r$ is odd. Then
\[ \begin{cases}
      x_{r}^{i}=\sum_{k=1}^{r-i}\gamma_{r,i+k}(-1)^{k}+1 & \mbox{if $i$ is even,} \\
		  x_{r}^{i}=\sum_{k=1}^{r-i}\gamma_{r,i+k}(-1)^{k+1}+1 & \mbox{if $i$ is odd.}
    \end{cases}
\]
\end{Cor}

\begin{Cor}
\label{correc2}
If $r$ is even. Then
\[ \begin{cases}
      x_{r}^{i}=\sum_{k=1}^{r-i}\gamma_{r,i+k}(-1)^{k} & \mbox{if $i$ is even,} \\
		  x_{r}^{i}=\sum_{k=1}^{r-i}\gamma_{r,i+k}(-1)^{k+1} & \mbox{if $i$ is odd.}
    \end{cases}
\]
\end{Cor}

Therefore, let $G$ be a cograph with cotree $T_{G_{r}}(a_{1},\ldots,a_{r-1},0|0,\ldots,0,a_{r})$. Then, by Corollaries \ref{correc1} and \ref{correc2}, we have that $-x_{r}^{i}$ is an eigenvalue of $G$ for $1\leq i\leq r-1$.


\section{Multiplicity}
\label{multi}

In this section, we compute the multiplicities of the eigenvalues $-x_{r}^{i}$ from Section \ref{recurrence} and we show how to retrieve the remaining eigenvalues of the cograph $G$ that has cotree $T_{G_{r}}(a_{1},\ldots,a_{r-1},0|0,\ldots,0,a_{r})$.

Suppose we have the cotree $T_{G_{i}}(a_{1},\ldots,a_{i-1},0|0,\ldots,0,a_{i})$.
If $i$ is even, every leave has value $x_{i}=0$, then we can apply \subiib~ of Algorithm Diagonalize $a_{1}a_{2}\cdots a_{i-1}(a_{i}-1)$ times. It creates  $a_{1}a_{2}\cdots a_{i-1}(a_{i}-1)$ permanent zeros in the removed leaves. And, the remaining leaves maintain the value zero as shown in Figure \ref{figobs1}. Therefore, $$m(-x_{r}^{i},G)\geq a_{1}a_{2}\cdots a_{i-1}(a_{i}-1).$$
\begin{figure}[h!]
\begin{center}
\begin{tikzpicture}
  [scale=1,auto=left,every node/.style={circle,scale=0.9}]

       \node[draw, circle, fill=blue!10, inner sep=0,label=left:$a_{i-1}$,label=above:$T_{G_{i}}$] (a) at (1,1) {$\cup$};
 			 \node[draw,circle,fill=black, label=below:$0$,label=left:$a_{i}$] (b) at (0,0) {};
			\node[draw,circle,fill=black, label=right:$x_{i}$,label=below:$0$,] (c) at (2,0) {};

				\node[draw, circle, fill=blue!10, inner sep=0,label=left:$a_{i-1}$] (d) at (5,1) {$\cup$};
 			 \node[draw,circle,fill=black, label=below:$0$,label=left:$a_{i}$] (e) at (5,0) {};

    			\node[draw, circle, fill=blue!10, inner sep=0,label=left:$a_{i-2}$,label=above:$T_{G_{i-1}}$] (f) at (9,1) {$\otimes$};
 			 \node[draw,circle,fill=black, label=below:$0$,label=left:$a_{i-1}$] (g) at (8,0) {};
			\node[draw,circle,fill=black, label=right:$x_{i-1}$,label=below:$0$,] (h) at (10,0) {};
							
							\draw[dashed, very thick] (0.3,0) -- (1.7,0);
							\draw[dashed, very thick] (8.3,0) -- (9.7,0);
							\draw [->] (6,0.5) -- (7,0.5);
                            \draw [->] (3,0.5) -- (4,0.5);
  \path
          (a) edge node[right]{}(b)
          (a) edge node[right]{}(c)
							
		   (d) edge node[right]{}(e)
           (f) edge node[right]{}(g)
           (f) edge node[right]{}(h)
        ;
\end{tikzpicture}
\caption{$i$ even.}
       \label{figobs1}
\end{center}			
\end{figure}

If $i$ is odd, every leave has value $x_{i}=1$, then we can apply \subib~ of Algorithm Diagonalize $a_{1}a_{2}\cdots a_{i-1}(a_{i}-1)$ times. It creates  $a_{1}a_{2}\cdots a_{i-1}(a_{i}-1)$ permanent zeros in the removed leaves. And, the remaining leaves maintain the value $1$ as shown in Figure \ref{figobs2}. Therefore, $$m(-x_{r}^{i},G)\geq a_{1}a_{2}\cdots a_{i-1}(a_{i}-1).$$

\begin{figure}[h!]
\begin{center}

\begin{tikzpicture}
  [scale=1,auto=left,every node/.style={circle,scale=0.9}]

       \node[draw, circle, fill=blue!10, inner sep=0,label=left:$a_{i-1}$,label=above:$T_{G_{i}}$] (a) at (1,1) {$\otimes$};
 			 \node[draw,circle,fill=black, label=below:$1$,label=left:$a_{i}$] (b) at (0,0) {};
			\node[draw,circle,fill=black, label=right:$x_{i}$,label=below:$1$,] (c) at (2,0) {};

				\node[draw, circle, fill=blue!10, inner sep=0,label=left:$a_{i-1}$] (d) at (5,1) {$\otimes$};
 			 \node[draw,circle,fill=black, label=below:$1$,label=left:$a_{i}$] (e) at (5,0) {};

    			\node[draw, circle, fill=blue!10, inner sep=0,label=left:$a_{i-2}$,label=above:$T_{G_{i-1}}$] (f) at (9,1) {$\cup$};
 			 \node[draw,circle,fill=black, label=below:$1$,label=left:$a_{i-1}$] (g) at (8,0) {};
			\node[draw,circle,fill=black, label=right:$x_{i-1}$,label=below:$1$,] (h) at (10,0) {};
							
							\draw[dashed, very thick] (0.3,0) -- (1.7,0);
							\draw[dashed, very thick] (8.3,0) -- (9.7,0);
							\draw [->] (6,0.5) -- (7,0.5);
                            \draw [->] (3,0.5) -- (4,0.5);
  \path
          (a) edge node[right]{}(b)
          (a) edge node[right]{}(c)
							
		   (d) edge node[right]{}(e)
           (f) edge node[right]{}(g)
           (f) edge node[right]{}(h)
        ;
\end{tikzpicture}
\caption{$i$ odd.}
       \label{figobs2}
\end{center}
\end{figure}

\begin{Lem}
\label{lem01}
Given $T_{G}(a_{1},\ldots,a_{r-1},0|0,\ldots,0,a_{r})$ and its eigenvalue $-x_{r}^{i}$. Then after the level $i$ the Algorithm Diagonalize with input $(T_{G},x_{r}^{i})$ does not create permanent zeros in the leaves.
\end{Lem}
{\bf Proof:} Consider the cotree $T_{G_{i}}(a_{1},\ldots,a_{i-1},0|0,\ldots,0,a_{i})$.
If $i$ is even, then $x_{i}=0$ and the cotree in the next level is $T_{G_{i-1}}(a_{1},\ldots,a_{i-2},0|0,\ldots,0,a_{i-1})$ with each leave $x_{i-1}=0$ as in Figure \ref{figlem1}.
\begin{figure}[h!]
\begin{center}
\begin{tikzpicture}
  [scale=1,auto=left,every node/.style={circle,scale=0.9}]

    		\node[draw, circle, fill=blue!10, inner sep=0,label=left:$a_{i-2}$,label=above:$T_{G_{i-1}}$] (f) at (9,1) {$\otimes$};
 			\node[draw,circle,fill=black, label=below:$0$,label=left:$a_{i-1}$] (g) at (8,0) {};
			\node[draw,circle,fill=black, label=right:$x_{i-1}$,label=below:$0$,] (h) at (10,0) {};
							
								\draw[dashed, very thick] (8.3,0) -- (9.7,0);
							
  \path

           (f) edge node[right]{}(g)
           (f) edge node[right]{}(h)
        ;
\end{tikzpicture}
\caption{$i$ even.}
       \label{figlem1}

\end{center}
\end{figure}

Using Theorem \ref{main2}, we have $m(0,G_{i-1})=0$, where $G_{i-1}$ is the cograph that has balanced cotree $T_{G_{i-1}}(a_{1},\ldots,a_{i-2},0|0,\ldots,0,a_{i-1})$. It means that, in the next iterations the permanent values in the leaves are different of zero.

If $i$ is odd, then $x_{i}=1$ and the cotree in the next level is $T_{G_{i-1}}$ $(a_{1},\ldots,a_{i-2},0|$ $0,\ldots,0,a_{i-1})$ with each leave $x_{i-1}=1$ as in Figure \ref{figlem2}.
\begin{figure}[h!]
\begin{center}
\begin{tikzpicture}
  [scale=1,auto=left,every node/.style={circle,scale=0.9}]

    		\node[draw, circle, fill=blue!10, inner sep=0,label=left:$a_{i-2}$,label=above:$T_{G_{i-1}}$] (f) at (9,1) {$\cup$};
 			\node[draw,circle,fill=black, label=below:$1$,label=left:$a_{i-1}$] (g) at (8,0) {};
			\node[draw,circle,fill=black, label=right:$x_{i-1}$,label=below:$1$,] (h) at (10,0) {};
							
								\draw[dashed, very thick] (8.3,0) -- (9.7,0);
							
  \path

           (f) edge node[right]{}(g)
           (f) edge node[right]{}(h)
        ;
\end{tikzpicture}
\caption{$i$ odd.}
       \label{figlem2}
\end{center}
\end{figure}
Using Theorem \ref{main2}, we have that $m(-1,G_{i-1})=0$, where $G_{i-1}$ is the cograph that has balanced cotree $T_{G_{i-1}}(a_{1},\ldots,a_{i-2},0|0,\ldots,0,a_{i-1})$. It means that, in the next iterations the permanent values in the leaves are different of zero.

Therefore, Algorithm Diagonalize with input $(T_{G_{r}},x_{r}^{i})$ does not create permanent zeros from $T_{G_{i-1}}(a_{1},\ldots,a_{i-2},0|0,\ldots,0,a_{i-1})$ onwards.

\begin{Thr}
For $1\leq i\leq r-1$ we have $$m(-x_{r}^{i},G)=a_{1}\cdots(a_{i}-1).$$
\end{Thr}
{\bf Proof:} Suppose we are at the cotree $T_{G_{i}}$ with leaves $x_{i}$.
According Lemma \ref{lema00}, in the previous steps the removed leaves have non-null permanent values. And, after the iterations at level $i$, Lemma \ref{lem01} guarantee that no new permanent zeros are created from the cotree $T_{G_{i-1}}$ onwards. Therefore, $m(-x_{r}^{i},G)=a_{1}\cdots(a_{i}-1)$.

Let $G$ be a cograph with balanced cotree $T_{G_{r}}(a_{1},\ldots,a_{r-1},0|0,\ldots,a_{r})$ of order $n=a_{1}a_{2}\cdots a_{r}$. By Lemma \ref{main2}, we know if $r$ is odd then $m(-1,G)=a_{1}a_{2}\cdots a_{r-1}(a_{r}-1)$ and, if $r$ is even then $m(0,G)=a_{1}a_{2}\cdots a_{r-1}(a_{r}-1)$. In both cases, odd or even, we have

$(a_{1}-1)+a_{1}(a_{2}-1)+a_{1}a_{2}(a_{3}-1)+\cdots+a_{1}\cdots a_{r-2}(a_{r-1}+1)+a_{1}\cdots a_{r-1}(a_{r}-1)=n-1$

integral engenvalues of the cograph. Then, the remaining eigenvalue can be computed using the trace of the adjacency matrix.

\begin{Cor}
\label{eigenvalue}
The remaining eigenvalue is
\[
\lambda=\begin{cases}
    -\sum_{i=1}^{n}m(x_{r}^{i})x_{r}^{i}              & \mbox{ if $r$ is even,} \\
    -\sum_{i=1}^{n}m(x_{r}^{i})x_{r}^{i}-(-1)(a_{r}-1)& \mbox{ if $r$ is odd.}
\end{cases}
\]

\end{Cor}
{\bf Proof:} Just use the fact that $Tr(A)=0$.


\section{Examples}
\label{examplesection}

In this section, we illustrate our results with two examples.

\begin{Ex}
\label{example1}
Let $G$ be a  cograph with balanced cotree $T_G  (2,2,2,0| 0,0,0,3)$ of order $n =2\cdot 2\cdot 2\cdot 3=24$ and $r=4$.
\begin{figure}[h!]
\begin{center}
\begin{tikzpicture}
  [scale=1,auto=left,every node/.style={circle,scale=0.9}]

                \node[draw, circle, fill=blue!10, inner sep=0] (a) at (5.75,4) {$\otimes$};
								
              \node[draw, circle, fill=blue!10, inner sep=0] (b) at (2.75,3) {$\cup$};
							 \node[draw, circle, fill=blue!10, inner sep=0] (c) at (8.75,3) {$\cup$};

        \node[draw, circle, fill=blue!10, inner sep=0] (d) at (1.25,2) {$\otimes$};
				\node[draw, circle, fill=blue!10, inner sep=0] (e) at (4.25,2) {$\otimes$};
				\node[draw, circle, fill=blue!10, inner sep=0] (f) at (7.25,2) {$\otimes$};
				\node[draw, circle, fill=blue!10, inner sep=0] (g) at (10.25,2) {$\otimes$};

       \node[draw, circle, fill=blue!10, inner sep=0] (h) at (0.5,1) {$\cup$};
       \node[draw, circle, fill=blue!10, inner sep=0] (i) at (2,1) {$\cup$};			
			\node[draw, circle, fill=blue!10, inner sep=0] (j) at (3.5,1) {$\cup$};
			\node[draw, circle, fill=blue!10, inner sep=0] (k) at (5,1) {$\cup$};
			\node[draw, circle, fill=blue!10, inner sep=0] (l) at (6.5,1) {$\cup$};
			\node[draw, circle, fill=blue!10, inner sep=0] (m) at (8,1) {$\cup$};
			\node[draw, circle, fill=blue!10, inner sep=0] (n) at (9.5,1) {$\cup$};
			\node[draw, circle, fill=blue!10, inner sep=0] (o) at (11,1) {$\cup$};
 			
			\node[draw,circle,fill=black] (h1) at (0,0) {};
			\node[draw,circle,fill=black] (h2) at (0.5,0) {};
			\node[draw,circle,fill=black] (h3) at (1,0) {};
			
			\node[draw,circle,fill=black] (i1) at (1.5,0) {};
			\node[draw,circle,fill=black] (i2) at (2,0) {};
			\node[draw,circle,fill=black] (i3) at (2.5,0) {};
			
			\node[draw,circle,fill=black] (j1) at (3,0) {};
			\node[draw,circle,fill=black] (j2) at (3.5,0) {};
			\node[draw,circle,fill=black] (j3) at (4,0) {};
			
			\node[draw,circle,fill=black] (k1) at (4.5,0) {};
			\node[draw,circle,fill=black] (k2) at (5,0) {};
			\node[draw,circle,fill=black] (k3) at (5.5,0) {};
			
			\node[draw,circle,fill=black] (l1) at (6,0) {};
			\node[draw,circle,fill=black] (l2) at (6.5,0) {};
			\node[draw,circle,fill=black] (l3) at (7,0) {};
			
			\node[draw,circle,fill=black] (m1) at (7.5,0) {};
			\node[draw,circle,fill=black] (m2) at (8,0) {};
			\node[draw,circle,fill=black] (m3) at (8.5,0) {};
			
			\node[draw,circle,fill=black] (n1) at (9,0) {};
			\node[draw,circle,fill=black] (n2) at (9.5,0) {};
			\node[draw,circle,fill=black] (n3) at (10,0) {};
			
			\node[draw,circle,fill=black] (o1) at (10.5,0) {};
			\node[draw,circle,fill=black] (o2) at (11,0) {};
			\node[draw,circle,fill=black] (o3) at (11.5,0) {};

			  \path
				  (a) edge node[right]{}(b)
					(a) edge node[right]{}(c)
					
					(b) edge node[right]{}(d)
					(b) edge node[right]{}(e)
					
					(c) edge node[right]{}(f)
					(c) edge node[right]{}(g)
				
				  (d) edge node[right]{}(h)
					(d) edge node[right]{}(i)
					
					(e) edge node[right]{}(j)
					(e) edge node[right]{}(k)
					
					(f) edge node[right]{}(l)
					(f) edge node[right]{}(m)
					
					(g) edge node[right]{}(n)
					(g) edge node[right]{}(o)
				
				  (h) edge node[right]{}(h1)
					(h) edge node[right]{}(h2)
					(h) edge node[right]{}(h3)
					
          (i) edge node[right]{}(i1)
					(i) edge node[right]{}(i2)
					(i) edge node[right]{}(i3)
					
					(j) edge node[right]{}(j1)
					(j) edge node[right]{}(j2)
					(j) edge node[right]{}(j3)
					
					(k) edge node[right]{}(k1)
					(k) edge node[right]{}(k2)
					(k) edge node[right]{}(k3)
					
					(l) edge node[right]{}(l1)
					(l) edge node[right]{}(l2)
					(l) edge node[right]{}(l3)
					
					(m) edge node[right]{}(m1)
					(m) edge node[right]{}(m2)
					(m) edge node[right]{}(m3)
					
					(n) edge node[right]{}(n1)
					(n) edge node[right]{}(n2)
					(n) edge node[right]{}(n3)
					
					(o) edge node[right]{}(o1)
					(o) edge node[right]{}(o2)
					(o) edge node[right]{}(o3)
							
        ;
\end{tikzpicture}
\caption{$T_{G_{4}}(2,2,2,0| 0,0,0,3)$.}
       \label{ex_r_even}
\end{center}			
\end{figure}
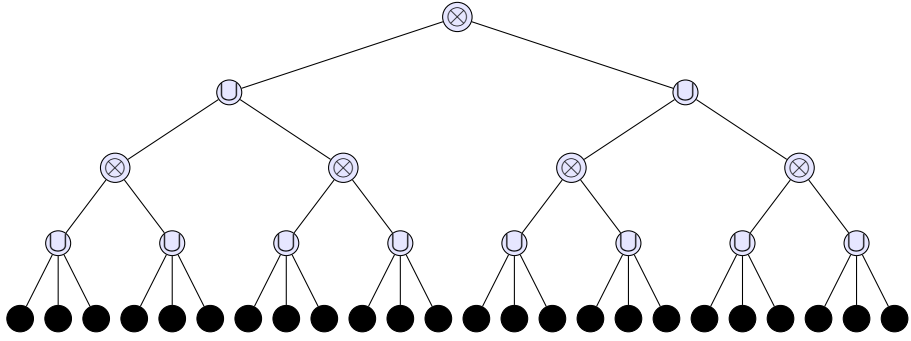	

Using Corollary \ref{correc2} in the cotree in Figure \ref{ex_r_even}. At each level $i$  for $1\leq i\leq 3$, we obtain the following.	

Level $i=1$: $$x_{4}^{1}=\sum_{k=1}^{3}\gamma_{4,1+k}(-1)^{k+1}=\gamma_{4,2}(-1)^{2}+\gamma_{4,3}(-1)^{3}+\gamma_{4,4}(-1)^{4}=a_{4}a_{3}a_{2}-a_{4}a_{3}+a_{4}=9.$$
Then the eigenvalue is $\lambda=-x_{4}^{1}=-9$. And its multiplicity is $m(\lambda)=a_{1}-1=1$.

Level $i=2$: $$x_{4}^{2}=\sum_{k=1}^{2}\gamma_{4,2+k}(-1)^{k}=\gamma_{4,3}(-1)^{1}+\gamma_{4,4}(-1)^{4}=-a_{4}a_{3}+a_{4}=-3.$$
Then the eigenvalue is $\lambda=-x_{4}^{2}=3$. And its multiplicity is $m(\lambda)=a_{1}(a_{2}-1)=2$.

Level $i=3$:
$$x_{4}^{3}=\sum_{k=1}^{1}\gamma_{4,3+k}(-1)^{k+1}=\gamma_{4,4}(-1)^{2}=a_{4}=3.$$
Then the eigenvalue is $\lambda=-x_{4}^{3}=-3$. And its multiplicity is $m(\lambda)=a_{1}a_{2}(a_{3}-1)=4$.

Using Theorem \ref{main2} we obtain $m(0,G)=16$.
Then, we still have to compute one eigenvalue. As in Corollary \ref{eigenvalue}, we use the trace of the adjacency matrix as follows.
 $$Tr(A)=0=0(16)+(-3)(4)+(3)(2)+(-9)(1)+\lambda$$
and it implies that $\lambda=15$. Therefore,
$$\mbox{Spec}(G)=\{(-9)^{(1)},(-3)^{(4)},(0)^{(16)},(3)^{(2)},(15)^{(1)}\}.$$

\end{Ex}

\begin{Ex}
\label{example2}
Let $G$ be a  cograph with balanced cotree $T_{G_{3}}(2,2,0| 0,0,2)$ of order $n =2\cdot 2\cdot 2=8$ and $r=3$.
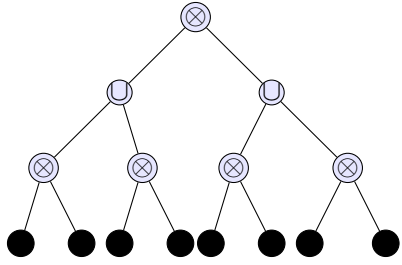
\begin{figure}[h!]
\begin{center}
\begin{tikzpicture}
  [scale=1,auto=left,every node/.style={circle,scale=0.9}]

  \node[draw,circle,fill=black,label=below:$$] (o) at (3,4) {};
  \node[draw,circle,fill=black,label=below:$$] (n) at (1.8,4) {};
   \node[draw,circle,fill=black,label=below:$$] (q) at (2.2,4) {};

  \node[draw, circle, fill=blue!10, inner sep=0] (m) at (2.5,5) {$\otimes$};
  \node[draw,circle,fill=blue!10, inner sep=0] (l) at (4,5) {$\otimes$};

  \node[draw,circle,fill=black,label=below:$$] (p) at (3.5,4) {};
  \node[draw,circle,fill=black,label=below:$$] (u) at (4.5,4) {};

  \node[draw, circle, fill=blue!10, inner sep=0] (j) at (3,6) {$\cup$};
  \node[draw,circle,fill=blue!10, inner sep=0] (h) at (2,7) {$\otimes$};
  \node[draw, circle, fill=blue!10, inner sep=0] (g) at (1,6) {$\cup$};
  \node[draw,circle,fill=blue!10, inner sep=0] (f) at  (1.3,5) {$\otimes$};
  \node[draw,circle,fill=blue!10, inner sep=0] (a) at (0,5) {$\otimes$};
  \node[draw, circle, fill=black, label=below:$$] (b) at (-0.3,4) {};
  \node[draw,circle,fill=black, label=below:$$] (c) at (0.5,4) {};
  \node[draw,circle,fill=black,label=below:$$] (e) at (1,4) {};

  \path (a) edge node[left]{} (b)
        (a) edge node[below]{} (c)

   (f) edge node[below]{} (e)
  (f) edge node[below]{} (n)

    (l) edge node[below]{} (p)
  (l) edge node[below]{} (u)

        (f) edge node[right]{}(g)
        (g) edge node[left]{}(a)
        (h) edge node[right]{}(j)
        (h) edge node[left]{}(g)

        (j) edge node[right]{}(l)
        (j) edge node[below]{}(m)
        (m) edge node[right]{}(o)
        (m) edge node[left]{} (q);
\end{tikzpicture}
       \caption{$T_{G_{3}}(2,2,0|0,0,2)$.}
       \label{Ex_2}
 \end{center}      
 \end{figure}

Using Corollary \ref{correc1} in the cotree in Figure \ref{Ex_2}. At each level $i$  for $1\leq i\leq 2$, we obtain the following.

Level $i=1$: $$x_{3}^{1}=\sum_{k=1}^{2}\gamma_{3,1+k}(-1)^{k+1}+1=\gamma_{3,2}(-1)^{2}+\gamma_{3,3}(-1)^{3}+1=a_{3}a_{2}-a_{3}
+1=3.$$
Then the eigenvalue is $\lambda=-x_{3}^{1}=-3$. And its multiplicity is $m(\lambda)=a_{1}-1=1$.

Level $i=2$: $$x_{3}^{2}=\sum_{k=1}^{1}\gamma_{3,2+k}(-1)^{k}+1=\gamma_{3,3}(-1)^{1}+1=-a_{3}+1=-1.$$
Then the eigenvalue is $\lambda=-x_{3}^{2}=1$. And its multiplicity is $m(\lambda)=a_{1}(a_{2}-1)=2$.

Using Theorem \ref{main2} we obtain $m(-1,G)=4$. Then, we still have to compute one eigenvalue. As in Corollary \ref{eigenvalue}, we use the trace of the adjacency matrix as follows.
$$Tr(A)=0=-3+1(2)-1(4)+\lambda$$
and it implies that $\lambda=5$. Therefore,
$$\mbox{Spec}(G)=\{(-3)^{(1)},(-1)^{(4)},(1)^{(2)},(5)^{(1)}\}.$$

\end{Ex}

\section{Application}
\label{applic}

\label{inter}
In this section, we show how to estimate the eigenvalues  of any cograph using a balanced cotree  $T_{G_{r}}(a_{1},\ldots,a_{r-1},0|0,\ldots,0,a_{r})$.

The following known result is an eigenvalue interlacing theorem and can be found in \cite{hall}.

\begin{Thr}
\label{interlacing}
Suppose $A \in \mathbb{R}^{n \times n}$ is symmetric. Let $B \in \mathbb{R}^{m\times m}$ with $m< n$ be a principal  submatrix (obtained by deleting both $i$-th row and
$i$-th column for some values of $i$). Suppose $A$ has eigenvalues $\lambda_1 \geq \ldots \geq \lambda_n$ and $B$ has eigenvalues $\beta_1 \geq \ldots \geq \beta_m$. Then
$$ \lambda_k  \geq \beta_k \geq \lambda_{k+n-m}   \hspace{0,75cm} for \hspace{0,3cm} k = 1, \ldots, m.$$
And if $m=n-1$,
$$ \lambda_1 \geq \beta_1 \geq \lambda_2 \geq \beta_2 \geq \ldots \geq \beta_{n-1} \geq \lambda_n$$
\end{Thr}

The balanced cotree that we use to estimate the eigenvalues of the cograph $G$ is obtained from the original cotree $T_{G}$ with depth $r$ by adding superfluous intermediate nodes so that  all leaves are in the same level. If a leaf $v$ is not at the level $r$, say it is at the level $s<r$, then we replace the edge hanging $v$ by a path of length $r-s$ with $v$ hanging, in such a way that the intermediate nodes alternate between $\cup$ and  $\otimes$. We finalize this process by adding (if necessary) new  superfluous intermediate nodes and duplicate or coduplicate vertices until that the cotree $T_{G}$ becomes into a balanced cotree $T_{G_{r}}$ associated to $T_{G}$.

\begin{Ex}
The standard cotree $T_{G}$ in Figure \ref{fig8} is associated to the cograph of Figure \ref{cotree}. And, the balanced cotree $T_{G_{3}}(2,2,0|0,0,2)$ associated to $T_{G}$ is shown in Figure \ref{fig9}.

\begin{figure}[h!]
       \begin{minipage}[c]{0.45 \linewidth}
\begin{tikzpicture}
 [scale=1,auto=left,every node/.style={circle,scale=0.9}]

  \node[draw,circle,fill=black,label=below:$$] (n) at (1.8,4) {};

  \node[draw, circle, fill=black,label=below:$$] (m) at (2.5,5) {};
  \node[draw,circle,fill=black,label=below:$$] (l) at (4,5) {};


  \node[draw, circle, fill=blue!10, inner sep=0] (j) at (3,6) {$\cup$};
  \node[draw,circle,fill=blue!10, inner sep=0] (h) at (2,7) {$\otimes$};
  \node[draw, circle, fill=blue!10, inner sep=0] (g) at (1,6) {$\cup$};
  \node[draw,circle,fill=blue!10, inner sep=0] (f) at  (1.3,5) {$\otimes$};
  \node[draw,circle,fill=blue!10, inner sep=0] (a) at (0,5) {$\otimes$};
  \node[draw, circle, fill=black, label=below:$$] (b) at (-0.3,4) {};
  \node[draw,circle,fill=black, label=below:$$] (c) at (0.5,4) {};
  \node[draw,circle,fill=black,label=below:$$] (e) at (1,4) {};

  \path (a) edge node[left]{} (b)
        (a) edge node[below]{} (c)

   (f) edge node[below]{} (e)
  (f) edge node[below]{} (n)


        (f) edge node[right]{}(g)
        (g) edge node[left]{}(a)
        (h) edge node[right]{}(j)
        (h) edge node[left]{}(g)

        (j) edge node[right]{}(l)
        (j) edge node[below]{}(m);
\end{tikzpicture}
\caption{$T_{G}$.}
       \label{fig8}
       \end{minipage}\hfill
       \begin{minipage}[c]{0.45 \linewidth}
     \begin{tikzpicture}
  [scale=1,auto=left,every node/.style={circle,scale=0.9}]

  \node[draw,circle,fill=black,label=below:$$] (o) at (3,4) {};
  \node[draw,circle,fill=black,label=below:$$] (n) at (1.8,4) {};
   \node[draw,circle,fill=black,label=below:$$] (q) at (2.2,4) {};

  \node[draw, circle, fill=blue!10, inner sep=0] (m) at (2.5,5) {$\otimes$};
  \node[draw,circle,fill=blue!10, inner sep=0] (l) at (4,5) {$\otimes$};

  \node[draw,circle,fill=black,label=below:$$] (p) at (3.5,4) {};
  \node[draw,circle,fill=black,label=below:$$] (u) at (4.5,4) {};

  \node[draw, circle, fill=blue!10, inner sep=0] (j) at (3,6) {$\cup$};
  \node[draw,circle,fill=blue!10, inner sep=0] (h) at (2,7) {$\otimes$};
  \node[draw, circle, fill=blue!10, inner sep=0] (g) at (1,6) {$\cup$};
  \node[draw,circle,fill=blue!10, inner sep=0] (f) at  (1.3,5) {$\otimes$};
  \node[draw,circle,fill=blue!10, inner sep=0] (a) at (0,5) {$\otimes$};
  \node[draw, circle, fill=black, label=below:$$] (b) at (-0.3,4) {};
  \node[draw,circle,fill=black, label=below:$$] (c) at (0.5,4) {};
  \node[draw,circle,fill=black,label=below:$$] (e) at (1,4) {};

  \path (a) edge node[left]{} (b)
        (a) edge node[below]{} (c)

   (f) edge node[below]{} (e)
  (f) edge node[below]{} (n)

    (l) edge node[below]{} (p)
  (l) edge node[below]{} (u)

        (f) edge node[right]{}(g)
        (g) edge node[left]{}(a)
        (h) edge node[right]{}(j)
        (h) edge node[left]{}(g)

        (j) edge node[right]{}(l)
        (j) edge node[below]{}(m)
        (m) edge node[right]{}(o)
        (m) edge node[left]{} (q);
\end{tikzpicture}
       \caption{$T_{G_{3}}(2,2,0|0,0,2)$.}
       \label{fig9}
        \end{minipage}
\end{figure}

   By Example \ref{example2}, we know $Spec(G_{3})=\{(-1)^{(4)},(1^{(2)}),(3)^{(1)},(5)^{(1)}\}$. And, using  Theorem \ref{main2}, we obtain $m(-1,G)=2$ and $m(0,G)=1$.
Theorem \ref{interlacing} leads to the following set of inequations.

\begin{eqnarray}
\label{int1}
  5  =  \lambda_{1}  \geq & \beta_{1} & \geq  \lambda_{3}  =  1 \\
\label{int2}
  3  =  \lambda_{2}  \geq & \beta_{2} & \geq  \lambda_{4}  =  1 \\
\label{int3}
  1  =  \lambda_{3}  \geq & \beta_{3} & \geq  \lambda_{5}  =  -1 \\
\label{int4}
  1  =  \lambda_{4}  \geq & \beta_{4} & \geq  \lambda_{6}  =  -1 \\
\label{int5}
  -1  =  \lambda_{5}  \geq & \beta_{5} & \geq \lambda_{7}  =  -1 \\
\label{int6}
  -1  =  \lambda_{6}  \geq & \beta_{6} & \geq \lambda_{8}  =  -1.
\end{eqnarray}

Inequations (\ref{int5}) and (\ref{int6}) imply that $\beta_{5}=\beta_{6}=-1$. Using the facts that no cograph has eigenvalue in the interval $(-1,0)$ and $m(0,G)=1$, inequations (\ref{int3}) and (\ref{int4}) lead to $\beta_{3}=0$ and $0<\beta_{4}\leq 1$. And, inequations (\ref{int1}) and (\ref{int2}) give  $1\leq\beta_{1}\leq 5$ and $1\leq\beta_{2}\leq 3$.
\end{Ex}

By construction, the following observation can be stated.

Let $T_G$ be the standard cotree of depth $r$ of a cograph $G$.
Then a balanced cotree $T_{G_{r}}(a_1, a_2, \ldots, a_{r-1},0|0,\ldots, a_r)$ can be obtained from $T_G$ by adding

\begin{enumerate}
\item[(1)] new superfluous intermediate nodes,
\item[(2)] duplicate or coduplicate  vertices into the terminal vertices.
\end{enumerate}

\begin{Ex}

Let $G$ be the cograph given by $G= (( t k_1) \otimes (t k_1)) \otimes k_{p}$ with $t\geq 2$, $p\geq 1$ and $t>p$. $T_G$ is the standard cotree of the cograph $G$, Figure \ref{fig10}, and its  balanced cotree associated $T_{G_{2}}( p+2,0|0, t)$ is shown  in Figure \ref{fig11}.

\begin{figure}[h!]
       \begin{minipage}[c]{0.45 \linewidth}
\begin{tikzpicture}
 [scale=1,auto=left,every node/.style={circle,scale=0.9}]

  \node[draw, circle, fill=black,label=below:$1$] (m) at (2,5) {};
  \node[draw,circle,fill=black,label=below:$t$] (l) at (3,5) {};


  \node[draw, circle, fill=blue!10, inner sep=0] (j) at (2.5,6) {$\cup$};

 \node[draw, circle, fill=black,label=below:$1$] (c) at (3.5,6) {};
\node[draw, circle, fill=black,label=below:$p$] (d) at (4.5,6) {};

  \node[draw,circle,fill=blue!10, inner sep=0] (h) at (2,7) {$\otimes$};
  \node[draw, circle, fill=blue!10, inner sep=0] (g) at (1,6) {$\cup$};
  \node[draw,circle,fill=black,label=below:$t$] (f) at  (1.5,5) {};
  \node[draw,circle,fill=black,label=below:$1$] (a) at (0.5,5) {};

                \draw[dotted, very thick] (0.7,5) -- (1.3,5);
                \draw[dotted, very thick] (2.2,5) -- (2.8,5);
                \draw[dotted, very thick] (3.7,6) -- (4.3,6);

  \path
   (h) edge node[left]{}(c)
 (h) edge node[left]{}(d)

        (f) edge node[right]{}(g)
        (g) edge node[left]{}(a)
        (h) edge node[right]{}(j)
        (h) edge node[left]{}(g)

        (j) edge node[right]{}(l)
        (j) edge node[below]{}(m);

\end{tikzpicture}
\caption{$T_{G}$.}
       \label{fig10}
       \end{minipage}\hfill
       \begin{minipage}[c]{0.45 \linewidth}
     \begin{tikzpicture}
 [scale=1,auto=left,every node/.style={circle,scale=0.9}]


  \node[draw, circle, fill=black,label=below:$1$] (m) at (1.8,5) {};
  \node[draw,circle,fill=black,label=below:$t$] (l) at (2.6,5) {};

  \node[draw,circle,fill=black,label=below:$1$] (p) at (3.125,5) {};
 \node[draw,circle,fill=black,label=below:$t$] (p1) at (3.75,5) {};

  \node[draw,circle,fill=black,label=below:$1$] (q) at (4.125,5) {};
  \node[draw,circle,fill=black,label=below:$t$] (q1) at (4.75,5) {};

  \node[draw, circle, fill=blue!10, inner sep=0] (j) at (2,6) {$\cup$};

 \node[draw, circle, fill=blue!10, inner sep=0,label=left:$1$] (c) at (3,6) {$\cup$};
\node[draw, circle, fill=blue!10, inner sep=0,label=right:$p$] (d) at (4,6) {$\cup$};

  \node[draw,circle,fill=blue!10, inner sep=0] (h) at (2,7) {$\otimes$};
  \node[draw, circle, fill=blue!10, inner sep=0] (g) at (1,6) {$\cup$};
  \node[draw,circle,fill=black,label=below:$t$] (f) at  (1.3,5) {};
  \node[draw,circle,fill=black,label=below:$1$] (a) at (0.5,5) {};

         \draw[dotted, very thick] (0.7,5) -- (1.1,5);
         \draw[dotted, very thick] (2,5) -- (2.4,5);
         \draw[dotted, very thick] (3.325,5) -- (3.55,5);
         \draw[dotted, very thick] (4.325,5) -- (4.55,5);
         \draw[dotted, very thick] (3.3,6) -- (3.7,6);

  \path 


 (d) edge node[below]{} (q)
(d) edge node[below]{} (q1)
    (c) edge node[below]{} (p)
  (c) edge node[below]{} (p1)

   (h) edge node[left]{}(c)
 (h) edge node[left]{}(d)

        (f) edge node[right]{}(g)
        (g) edge node[left]{}(a)
        (h) edge node[right]{}(j)
        (h) edge node[left]{}(g)

        (j) edge node[right]{}(l)
        (j) edge node[below]{}(m);

\end{tikzpicture}
       \caption{$T_{G_{2}}$.}
       \label{fig11}
        \end{minipage}
\end{figure}

Using Corollary \ref{correc2}, we can compute $$Spec(G_{2})=\{(t(p+1))^{(1)},(0)^{((p+2)(t-1))},(-t)^{(p+1)}\}.$$
By Theorem \ref{main2}, we obtain  $m(0,G) =(2t-2)$ and  $m(-1,G)= p-1$.
We consider $m=|V(G)|=2t+p$ and $n=|V(G_{2})|=(p+2)t$.
Notice that, $m(0,G)+m(-1,G)=2t-2+p-1=m-3$. It means that, we still have to compute three eigenvalues of $G$. Suppose $\beta_{1}\geq\beta_{2}\geq\beta_{3}$. Using the interlacing theorem, we obtain the following.

$$\lambda_{1}\geq\beta_{1},\mbox{ it means that }t(p+1)\geq\beta_{1}.$$

$\lambda_{2}\geq\beta_{2},\mbox{ that is, }0\geq\beta_{2}$.
But $m(0,G)=2t-2$ and no cograph has eigenvalue in the interval $(-1,0)$. Then, we can conclude $$-1>\beta_{2}.$$

And, $\beta_{3}\geq\lambda_{t(p+2)}$, it implies that $$\beta_{3}\geq\ -t.$$

\end{Ex}

\section{Acknowledgements}
F. Tura acknowledges the support of FAPERGS (Grant 17/2551-0000813-8).

\end{document}